\documentclass[11pt,a4paper,oneside]{article}%
\usepackage{amsfonts}
\usepackage{graphicx}
\usepackage{amsmath,amssymb}
\usepackage{color}
\usepackage{amssymb}
\usepackage{amsmath}
\usepackage[T1]{fontenc}
\usepackage{amsthm}
\usepackage[margin=1in]{geometry}
\usepackage{amssymb}
\usepackage{amsfonts}
\usepackage{graphicx}
\usepackage{amsmath}%
\providecommand{\U}[1]{\protect\rule{.1in}{.1in}}
\newtheorem{theorem}{Theorem}

\newtheorem{corollary}[theorem]{Corollary}

\newtheorem{lemma}[theorem]{Lemma}

\newtheorem{proposition}[theorem]{Proposition}
\newtheorem{remark}[theorem]{Remark}

\begin{document}

\begin{center}
{\Large On the universality of spectral limit for random matrices with
martingale differences entries}

\bigskip

F. Merlev\`{e}de, C. Peligrad\footnote{Supported in part by a Charles Phelps
Taft Memorial Fund grant.} and M. Peligrad\footnote{Supported in part by a
Charles Phelps Taft Memorial Fund grant, and the NSF grant DMS-1208237.}

\bigskip

\end{center}

Universit\'{e} Paris Est, LAMA (UMR 8050), UPEM, CNRS, UPEC.

Email: florence.merlevede@u-pem.fr \bigskip

Department of Mathematical Sciences, University of Cincinnati, PO Box 210025,
Cincinnati, Oh 45221-0025, USA.

Email:peligrc@ucmail.uc.edu, peligrm@ucmail.uc.edu

\bigskip


\bigskip

\begin{center}
\bigskip\textbf{Abstract}
\end{center}

For a class of symmetric random matrices whose entries are martingale
differences adapted to an increasing filtration, we prove that under a
Lindeberg-like condition, the empirical spectral distribution behaves
asymptotically similarly to a corresponding matrix with independent centered
Gaussian entries having the same variances. Under a slightly reinforced
condition, the approximation holds in the almost sure sense. We also point out
several sufficient regularity conditions imposed to the variance structure for
convergence to the semicircle law or the Marchenko-Pastur law and other
convergence results. In the stationary case we obtain a full extension from
the i.i.d. case to the martingale case of the convergence to the semicircle
law as well as to the Marchenko-Pastur one. Our results are well adapted to
study several examples including non linear ARCH($\infty$) random fields.

\section{Introduction}

\qquad Some of the most celebrated theorems concerning the limiting density of
empirical spectral measure for large random matrices are Wigner's (1958)
semicircle law and Marchenko-Pastur (1967) law for covariance matrix. The
results have been extended in various directions. In the non-i.i.d. case
Pastur (1973) showed that a Lindeberg-like condition is sufficient for the
convergence to the semicircle law (see also Girko \textit{et al.} (1994) and
Girko (2013)). It was shown that the Lindeberg's condition is also relevant
for convergence to the Marchenko-Pastur law (see Theorem 3.10 in
Bai-Silverstein, 2010). Recently, Tao and Vu (2010) obtained the circular law
as spectral limit for matrices with independent entries. All these results
assume the independence between the entries of the matrix. An important
feature of these results is that the empirical spectral measure converges in
distribution for almost all points in the sample space.

For dependent entries the situation is not so well understood. Chatterjee
(2006) treated exchangeable entries. Several authors considered the martingale
difference type entries. Steps in this direction are papers by G\"{o}tze and
Tikhomirov (2004, 2006) and G\"{o}tze \textit{et al.} (2012) who treat the
semicircle law, and papers by Adamczak (2011, 2013) and O'Rourke (2012) who
deal with the Marchenko-Pastur law. These works study the universality for the
empirical distribution function when the martingale difference property is
defined for an entry of the matrix conditioned by the "past" which is not an
ordered filtration, so the results cannot be applied to several martingale
random fields useful in statistical applications. Furthermore, the conditions
imposed in the stationary case lead to constant conditional variance, with
respect to the "past".

There are many time series in econometric theory that can be modeled by an
autoregressive process with martingale innovations which have nonconstant
conditional variance (heteroscedasticity). A basic diagnostic for knowing that
such a model is adequate is to look at the Wachter plot (i.e. to plot the
values of the ordered eigenvalues against the quantiles of the
Marchenko--Pastur law or Wigner law). Our paper provides a theoretical
justification of such a procedure. Therefore, with a view towards
applications, the main goal of our paper is to study the universality problem
for a more general class of martingale differences which are adapted to an
increasing filtration. We also impose a mild mixing condition that allows us
to go beyond the constant conditional variance imposed in the previous
studies, making possible to treat models that present heteroscedasticity. We
provide two types of results, one concerning convergence in probability, and
another concerning convergence in distribution of the empirical spectral
density for almost all points in the sample space, which we believe is the
first one of this type for martingale dependences. As corollaries we point out
convergence to the semicircle law, the Marchenko-Pastur law as well as other
limits for the limiting spectral density. For martingale differences which are
selected from a stationary random field we obtain, without any additional
conditions, a generalization of the empirical spectral theorems for i.i.d. We
point out several applications of our results to ARCH models and matrices
constructed from a triangular array of one dimensional martingales.

Our method consists in comparing the Stieltjes transform of the random matrix
with martingale like entries with the Stieltjes transform of a Gaussian matrix
with the same covariance structure, which has interest in itself. The proofs
are based on a blend of Lindeberg-like method, blocking techniques and
delicate maximal inequalities. The blocking is needed to overcome the
difficulties raised by selecting meaningful filtrations and mixing conditions
associated to random fields.

The paper is organized in the following way. In Section 2 we list the
approximations results, spectral limit theorems, and provide a discussion of
our conditions. Applications are included in Section 3. Section 4 is devoted
to the main proofs. Finally, in Section 5, we carry out the proofs of some
technical results which are important in themselves and also provide some
background material.

All along the paper, for positive numbers $a_{n}$ and $b_{n}$, the notation
$a_{n}\ll b_{n}$ means that for a positive constant $c$, we have $a_{n}\leq
c\, b_{n}$ for all $n$.

\section{ Results}

Let $(X_{\ell k})_{(\ell,k)\in\mathbb{Z}^{2}}$ be real-valued random variables
such that $\mathbb{E}(X_{\ell k})=0$ and $\mathbb{E}(X_{\ell k}^{2}%
)=\sigma_{\ell k}^{2}$, and let $(Y_{ij})_{(i,j)\in{\mathbb{N}}^{2}}$ be a
sequence of independent centered real-valued Gaussian r.v.'s with
${\mathbb{E}}(Y_{ij}^{2})=\sigma_{ij}^{2}$ which is in addition independent of
$(X_{\ell k})_{(\ell,k)\in\mathbb{Z}^{2}}.$ We shall assume that the variables
are defined on the same probability space $(\Omega,\mathcal{F},\mathbb{P}).$

We consider the symmetric $n\times n$ random matrix $\mathbf{X}_{n}$ such
that, for any $i$ and $j$ in $\{1,\dots,n\}$
\begin{align}
(\mathbf{X}_{n})_{ij}  &  =X_{ij}\,\text{ for }i\geq j\ \text{ and
}\label{defX}\\
(\mathbf{X}_{n})_{ij}  &  =X_{ji}\,\text{ for }i<j\,.\nonumber
\end{align}
Denote by $\lambda_{1}^{n}\leq\dots\leq\lambda_{n}^{n}$ the eigenvalues of
\begin{equation}
{\mathbb{X}}_{n}:=\frac{1}{n^{1/2}}\mathbf{X}_{n} \label{defW}%
\end{equation}
and define its distribution function by
\[
\ \mathbf{F}^{{\mathbb{X}}_{n}}(t)=\frac{1}{n}{\sum\limits_{1\leq k\leq n}%
}I(\lambda_{k}\leq t)\,,
\]
where $I(A)$ denotes the indicator of an event $A$.

Similarly we define $\mathbf{Y}_{n}$ and ${\mathbb{Y}}_{n}$ and $\mathbf{F}%
^{{\mathbb{Y}}_{n}}(t).$

The Levy distance between two distribution functions $F$ and $G$ is defined
by
\[
d(F,G) = \inf\{ \varepsilon>0 \ : \ F(x-\varepsilon)-\varepsilon\leq G(x)\leq
F(x+\varepsilon)+\varepsilon\} \, .
\]
It is well-known that a sequence of distribution functions $F_{n}(x)$
converges to a distribution function $F(x)$ at all continuity points $x$ of
$F$ if and only if $d(F_{n},G)\rightarrow0.$ We shall refer to this
convergence as weak convergence and denote $F_{n}\Rightarrow F$. In this paper
we are interested in two types of results.

1. Convergence in probability. There is a distribution function $\mathbf{F}$
such that for all positive $\epsilon$%
\begin{equation}
\lim_{n\rightarrow\infty}\mathbb{P}(d(\mathbf{F}^{{\mathbb{X}}_{n}}%
,\mathbf{F})>\epsilon)=0 \, . \label{convprob}%
\end{equation}
By abusing the language, for simplicity, we shall denote this type of
convergence $\mathbf{F}^{{\mathbb{X}}_{n}}\Rightarrow\mathbf{F}$ in probability.

2. Convergence almost sure. There is a distribution function $\mathbf{F}$ such
that%
\begin{equation}
\mathbb{P}(\lim_{n\rightarrow\infty}d(\mathbf{F}^{{\mathbb{X}}_{n}}%
,\mathbf{F})=0)=1. \label{conva.s}%
\end{equation}
In the sequel the last convergence will be denoted $\mathbf{F}^{{\mathbb{X}%
}_{n}}\Rightarrow\mathbf{F}$ a.s.

\bigskip

The Stieltjes transform of $\mathbf{F}^{{\mathbb{X}}_{n}}$ is given by
\begin{equation}
S^{{\mathbb{X}}_{n}}(z)=\int\frac{1}{x-z}d\mathbf{F}^{{\mathbb{X}}_{n}%
}(x)=\frac{1}{n}\mathrm{Tr}(n^{-1/2}\mathbf{X}_{n}-z{\mathbf{I}}_{n})^{-1}\,,
\label{StTr}%
\end{equation}
where $z=u+iv\in\mathbb{C}^{+}$ (the set of complex numbers with positive
imaginary part), and $\mathbf{{I}}_{n}$ is the identity matrix of order $n$.

In order to introduce the filtration we shall use lexicographic order on
${\mathbb{Z}}^{2}$: if $\mathbf{i}=(i_{1},i_{2})$ and $\mathbf{j}=(j_{1}%
,j_{2})$ are distinct elements of ${\mathbb{Z}}^{2}$ the notation
$\mathbf{j}\leq_{\text{lex}}\mathbf{i}$ means that either $i_{1}\leq j_{1}$ or
$i_{1}=j_{1}$ and $i_{2}\leq j_{2}$ and the notation $\mathbf{j}<_{\text{lex}%
}\mathbf{i}$ means that either $i_{1}<j_{1}$ or $i_{1}=j_{1}$ and $i_{2}%
<j_{2}$. For any non-negative integer $a$, we introduce now a set of indexes%
\begin{equation}
B_{ij}^{a}=\{(u,v)\in{\mathbb{Z}}^{2};\mathit{\,\max}(|u-i|,|v-j|)\geq
a,\,(u,v)\leq_{\text{lex}}(i,j)\} \label{defB}%
\end{equation}
and for $i\geq j$ the filtration
\begin{align}
{\mathcal{F}}_{ij}^{a}  &  =\sigma(X_{uv}:(u,v)\in B_{ij}^{a}\ \text{and}%
\ v\leq u)\,\text{ if }B_{ij}^{a}\neq\emptyset\,,\label{deffiltration}\\
{\mathcal{F}}_{ij}^{a}  &  =\{\emptyset,\Omega\}\,\text{ if }B_{ij}%
^{a}=\emptyset\text{ and }\nonumber\\
{\mathcal{F}}_{ji}^{a}  &  ={\mathcal{F}}_{ij}^{a}\,.\nonumber
\end{align}
Note that $X_{ij}$\textit{ }is adapted to ${\mathcal{F}}_{ij}^{0},$ which is
an increasing filtration in lexicographic order. Our first result compares the
distribution of the spectral density of a matrix of martingale difference with
the spectral density of a matrix with Gaussian independent entries, defined
above. Here and everywhere in the paper we use the standard notation $\Vert
X\Vert_{p}=(\mathbb{E}|X|^{p})^{1/p}$ (for $X$ a real or complex-valued random variable).

\begin{theorem}
\label{replacement} Assume that for all $1\leq j\leq i$,
\begin{equation}
\mathbb{E}(X_{ij}|{\mathcal{F}}_{ij}^{1})=0\text{ a.s.} \label{mart}%
\end{equation}
and that
\begin{equation}
\sup_{n}\frac{1}{n^{2}}\sum_{1\leq j\leq i\leq n}\sigma_{ij}^{2}<\infty\,.
\label{bound}%
\end{equation}
Assume in addition that
\begin{equation}
\lim_{a\rightarrow\infty}\limsup_{n\rightarrow\infty}\frac{1}{n^{2}}%
\sum_{1\leq j\leq i\leq n}\Vert\mathbb{E}(X_{ij}^{2}-\sigma_{ij}%
^{2}|{\mathcal{F}}_{ij}^{a})\Vert_{1}=0\,, \label{mix}%
\end{equation}
and for any $\varepsilon>0$,
\begin{equation}
\frac{1}{n^{2}}\sum_{1\leq j\leq i\leq n}\mathbb{E}(X_{ij}^{2}I(|X_{ij}%
|>\varepsilon n^{1/2}))\rightarrow0\,. \label{L}%
\end{equation}
Then, for all $z\in\mathbb{C}^{+}$,
\begin{equation}
S^{{\mathbb{X}}_{n}}(z)-S^{{\mathbb{Y}}_{n}}(z)\rightarrow0\text{ in
probability.} \label{conclusion1}%
\end{equation}

\end{theorem}

Under a slightly stronger moment condition we obtain an almost sure result.

\begin{theorem}
\label{almost sure}Assume condition (\ref{mart}) is satisfied. Assume also
that for some non-decreasing function $h(x)\geq1$ such that $x^{-1}h(x)$ is
non-increasing and ${\sum\nolimits_{n}}(nh(n))^{-1}<\infty$, there exists a
positive constant $C$ such that
\begin{equation}
\sup_{i,j}\mathbb{E}(X_{ij}^{2}h(|X_{ij}|))\leq C\,, \label{higher moment}%
\end{equation}
and the following condition holds
\begin{equation}
\lim_{a\rightarrow\infty}\limsup_{n\rightarrow\infty}\frac{1}{n^{2}}%
\sum_{1\leq j\leq i\leq n}|\mathbb{E}(X_{ij}^{2}-\sigma_{ij}^{2}|{\mathcal{F}%
}_{ij}^{a})|=0\text{ a.s.} \label{mix a s}%
\end{equation}
Then for all $z\in{\mathbb{C}}^{+}$
\begin{equation}
S^{{\mathbb{X}}_{n}}(z)-S^{{\mathbb{Y}}_{n}}(z)\rightarrow0\text{ a.s.}
\label{conv a.s.}%
\end{equation}

\end{theorem}

\medskip

The relevance of these two theorems is that they make possible to transport
the limit results from Gaussian random matrices to matrices with martingale
structure. It is well known that in order to establish the convergence of
empirical spectral distribution of a sequence of matrices, one needs only to
show the convergence of their Stieltjes transforms and the limiting spectral
distribution can be obtained from the limiting Stieltjes transform (see
Theorem B.9 in Bai-Silverstein (2010), or Corollary 1 in Geronimo and Hill
(2003), combined with arguments on page 38 in Bai-Silverstein (2010), based on
Vitali's convergence theorem).

With the notations in definitions (\ref{convprob}) and (\ref{conva.s}), let us
give two corollaries of the above theorems:

\begin{corollary}
\label{corConvproba}Assume that $(X_{ij})_{(i,j)\in{\mathbb{Z}}^{2}}$ is as in
Theorem \ref{replacement}. Furthermore, assume that,
\[
\mathbf{F}^{{\mathbb{Y}}_{n}}\ \Rightarrow\mathbf{F}\ \text{ in probability,}%
\]
where $\mathbf{F}$ is a nonrandom distribution function. Then,
\[
\mathbf{F}^{{\mathbb{X}}_{n}}\Rightarrow\mathbf{F}\ \text{ in probability.}%
\]

\end{corollary}

The following corollary is a direct consequence of Theorem \ref{almost sure}
and Theorem B.9 in Bai-Silverstein (2010).

\begin{corollary}
\label{corConvas} Assume that $(X_{ij})_{(i,j)\in{\mathbb{Z}}^{2}}$ is as in
Theorem \ref{almost sure}. Furthermore, assume that,
\[
\mathbf{F}^{{\mathbb{Y}}_{n}}\ \Rightarrow\mathbf{F}\ \text{ a.s.}%
\]
where $\mathbf{F}$ is a nonrandom distribution function. Then,
\[
\mathbf{F}^{{\mathbb{X}}_{n}}\ \Rightarrow\mathbf{F}\ \text{ a.s.}%
\]

\end{corollary}

\begin{remark}
\label{remark1} Our Theorem \ref{replacement} also holds if the random
variables $X_{ij}$ are replaced by a triangular array $X_{n,ij}$ with $\ j\leq
i\ .$ For this case the filtration is defined as ${\mathcal{F}}_{n,ij}%
^{a}=\sigma(X_{n,uv}:(u,v)\in B_{ij}^{a}\ $and$\ \ v\leq u\ ).$ The conditions
of Theorem \ref{replacement} should be modified accordingly, meaning that the
additional index $n$ should be added in all the conditions. \medskip
\end{remark}

\begin{remark}
By the contractivity properties of the conditional expectation, the conditions
in Theorem \ref{replacement} could be imposed to larger sigma algebras
${\mathcal{K}}_{ij}^{a}$ such that $\mathcal{F}_{ij}^{a}\subseteq{\mathcal{K}%
}_{ij}^{a}$. For the selection ${\mathcal{K}}_{ij}^{a}={\mathcal{F}}_{ij}^{1}$
for all $a$, condition (\ref{mix}) is implied by
\begin{equation}
\lim_{n\rightarrow\infty}\frac{1}{n^{2}}\sum_{1\leq j\leq i\leq n}%
\Vert\mathbb{E}(X_{ij}^{2}-\sigma_{ij}^{2}|{\mathcal{F}}_{ij}^{1})\Vert
_{1}=0\,, \label{limitation}%
\end{equation}
which is similar to G\"{o}tze \textit{et al.} (2012) martingale difference
condition but with a smaller filtration. The advantage of our condition
(\ref{mix}) is that is well adjusted to take care of martingale differences
which form a stationary random field. \medskip
\end{remark}

\begin{remark}
We cannot use the same simple argument to enlarge the filtration used in
Theorem \ref{almost sure}. However the proof of this theorem is based on
moment estimates and we notice that the conclusion of Theorem
\ref{almost sure} holds if we replace condition (\ref{mix a s}) by the
following condition:%
\[
\lim_{n\rightarrow\infty}\frac{1}{n^{2}}\sum_{1\leq j\leq i\leq n}%
|\mathbb{E}(X_{ij}^{2}-\sigma_{ij}^{2}|{\mathcal{F}}_{ij}^{1})|=0\ \text{a.s.}%
\]

\end{remark}

\begin{remark}
A careful analysis of the proof of Theorem \ref{almost sure} reveals that
under a stronger stationarity assumption, condition (\ref{higher moment}) can
be replaced by a weaker condition. More precisely, we infer that we can
replace condition (\ref{higher moment}) by the following one: There is a
random variable $X$ such that
\[
\sup_{i,j}{\mathbb{P}}(|X_{ij}|>x)\leq{\mathbb{P}}(|X|>x)\,,
\]
with
\[
\mathbb{E}(X^{2}\ln(1+|X|))<\infty\,.
\]
Furthermore, in the strictly stationary case we can assume only the existence
of moments of order two (see Theorem \ref{stationary}).\newline
\end{remark}

\smallskip

\noindent\textbf{Convergence results.} Our results can be combined with all
the available results for orthogonal Gaussian ensembles to obtain various
limiting laws.

\medskip

1. \textbf{Convergence to the semicircle law.}

Let $g(x)$ and $G(x)$ denote the density and the distribution function of the
standard semicircle law:%

\[
g(x)=\frac{1}{2\pi}\sqrt{4-x^{2}}I(|x|\leq2),\text{ }G(x)=\int_{-\infty}%
^{x}g(u)du.
\]
Combining Theorem \ref{replacement} with Theorem 1.1 in G\"{o}tze and
Tikhomirov (2004) we obtain under additional regularity condition the
following result:

\begin{corollary}
\label{semicircular} Assume besides the conditions of Theorem
\ref{replacement} that
\begin{equation}
\frac{1}{n^{2}}{\sum_{1\leq j\leq i\leq n}}|\sigma_{ij}^{2}-1|\rightarrow0\,.
\label{Bogus0}%
\end{equation}
Then,
\[
\mathbf{F}^{{\mathbb{X}}_{n}}\Rightarrow G\text{ in probability .}
\]

\end{corollary}

\begin{corollary}
\label{semicircular a.s.}If the conditions of Theorem \ref{almost sure} and
(\ref{Bogus0}) are satisfied then,
\[
\ \mathbf{F}^{{\mathbb{X}}_{n}}\Rightarrow G\text{ a.s. }%
\]

\end{corollary}

We consider next a symmetric random matrix which is constructed with variables
$(X_{ij})_{1\leq j\leq i\leq n}$ from a stationary real-valued random field
$(X_{\mathbf{u}})_{\mathbf{u}\in\mathbb{Z}^{2}}.$ This means that for all $n$
and any $\mathbf{t}, \mathbf{u}_{1}, \dots, \mathbf{u}_{n}$ in ${\mathbb{Z}%
}^{2}$ such that $\mathbf{u}_{1}<_{\text{lex}}\mathbf{u}_{2}<_{\text{lex}%
}...<_{\text{lex}}\mathbf{u}_{n}$, $(X_{\mathbf{u}_{1}},X_{\mathbf{u}_{2}%
},..,X_{\mathbf{u}_{n}})$ has the same distribution as $(X_{\mathbf{u}%
_{1}+\mathbf{t}},X_{\mathbf{u}_{2}+\mathbf{t}},..,X_{\mathbf{u}_{n}%
+\mathbf{t}})$.

In this case we have the following generalization of the semicircle law from
an i.i.d. to the martingale difference sequences:

\begin{theorem}
\label{stationary} Assume that ${\mathbb{X}}_{n}$ is defined by (\ref{defW})
and based on a stationary real-valued random field $(X_{\mathbf{u}%
})_{\mathbf{u}\in\mathbb{Z}^{2}}$. Let $\mathcal{F}_{\mathbf{0}}^{\infty}%
=\cap_{a\in\mathbb{N}}{\mathcal{F} }_{\mathbf{0}}^{a}$ where $\mathbf{0}%
=(0,0)$. Assume that
\[
\mathbb{E}X_{\mathbf{0}}^{2}=1,\text{ }\mathbb{E(}X_{\mathbf{0}}|{\mathcal{F}%
}_{\mathbf{0}}^{1})=0\text{ a.s. and }\mathbb{E(}X_{\mathbf{0}}^{2}%
|{\mathcal{F}}_{\mathbf{0}}^{\infty})=1\text{ a.s.}\label{L2cond}%
\]
Then,
\[
\mathbf{F}^{{\mathbb{X}}_{n}}\ \Rightarrow G\ \text{ a.s.}%
\]

\end{theorem}

2. \textbf{Convergence to the Marchenko-Pastur law.}

The sample covariance matrix is very important in multivariate statistical
inference. Suppose we have real matrices $\mathbf{X}=\mathbf{X}_{np}%
=(X_{ij})_{1\leq i\leq p,1\leq j\leq n}$. The sample covariance matrix is
simply defined as%
\[
\mathbf{A=}\frac{1}{n}\mathbf{X}\mathbf{X}^{T}\,,
\]
where $\mathbf{X}^{T}$ is the transpose matrix of $\mathbf{X}$. We shall
assume that $p/n\rightarrow y$ where $y \in(0,\infty)$. In the context of
independent entries with the same mean, variance $1$ and satisfying (\ref{L})
(where the sum extends over $1\leq i\leq p$ and $1\leq j\leq n$), the limiting
spectral distribution follows the standard Marchenko-Pastur law with the
density%
\[
\label{M-P}\tilde{g}_{y}(x)=\frac{1}{2\pi xy}\sqrt{(c-x)(x-b)}I(b\leq x\leq
c)
\]
and a point mass $1-1/y$ at the origin if $y>1,$ where $b=(1-\sqrt{y})^{2}$
and $c=(1+\sqrt{y})^{2}.$ See Theorem 3.10 in Bai-Silverstein (2010) and the
references therein.

It is well-known that for deriving the limiting spectral distribution of
$\mathbf{A}$ it is enough to study the Stieltjes transform of the following
symmetric matrix of order $N=n+p$:
\[
\mathbf{B}_{N}=\frac{1}{\sqrt{n}}\left(
\begin{array}
[c]{cc}%
\mathbf{0} & \mathbf{X}^{T}\\
\mathbf{X} & \mathbf{0}%
\end{array}
\right)  \,.
\]
Indeed the eigenvalues of $\mathbf{B}_{N}^{2}$ are the eigenvalues of
$n^{-1}\mathbf{X}^{T}\mathbf{X}$ together with the eigenvalues of
$n^{-1}\mathbf{X}\mathbf{X}^{T}$. Assuming that $p\leq n$ (otherwise exchange
the role of $\mathbf{X}$ and $\mathbf{X}^{T}$ everywhere), the following
relation holds: for any $z\in{\mathbb{C}}^{+}$
\begin{equation}
S_{\mathbf{A}}(z)=z^{-1/2}\frac{N}{2p}S_{\mathbf{B}_{N}}(z^{1/2})+\frac
{p-n}{2pz}\,. \label{rel1betweenS}%
\end{equation}
This relationship together with our results make it possible to formulate the
convergence to Marchenko-Pastur law for martingale difference entries. For
instance we can give the following result which follows easily by using
Theorem \ref{replacement} together with Remark \ref{remark1}, applied to the
matrix $\mathbf{B}_{N}:=n^{-1/2} (b_{i,j})_{1\leq j \leq i \leq N}$ where
$b_{i,j}=X_{i-n,j}{\mathbf{1}}_{i\geq n+1}{\mathbf{1}}_{1 \leq j \leq n}$.

\begin{theorem}
\label{Th-Ma-Pa} Suppose we have matrices $\mathbf{X}=(X_{ij})_{1\leq i\leq
p,1\leq j\leq n}$ of centered, square integrable real-valued r.v.'s with the
same variance equals to $1$ and $p/n\rightarrow y$ where $y \in(0,\infty)$.
Assume that for all $(i,j)$ such that $1\leq i\leq p$ and $1\leq j\leq n$,
\[
\mathbb{E}(X_{ij}|\sigma(X_{\mathbf{u}};\mathbf{u}<_{\text{lex}}(i,j))=0\text{
a.s.}%
\]
Assume in addition that
\[
\lim_{a\rightarrow\infty}\limsup_{n\rightarrow\infty}\frac{1}{n^{2}}\sum
_{i=1}^{p}\sum_{j=1}^{n}\Vert\mathbb{E}(X_{ij}^{2}-1|\sigma(X_{\mathbf{u}%
};\mathbf{u}\in B_{ij}^{a}))\Vert_{1}=0\,,
\]
where $B_{ij}^{a}$ is defined by \eqref{defB}, and for any $\varepsilon>0$,
\[
\lim_{n\rightarrow\infty}\frac{1}{n^{2}}\sum_{i=1}^{p}\sum_{j=1}^{n}%
\mathbb{E}(X_{ij}^{2}I(|X_{ij}|>\varepsilon n^{1/2}))=0\,.
\]
Then,
\[
\mathbf{F}^{\mathbf{XX}^{T}/n}\ \Rightarrow\tilde{G}_{y}\ \ \text{\ in
probability}\,\text{,}%
\]
where $\tilde{G}_{y}$ is the standard Marchenko-Pastur distribution function.
\end{theorem}

When the entries of the matrices $\mathbf{X}=(X_{ij})_{1\leq i\leq p,1\leq
j\leq n}$ come from a stationary random field, we can formulate an almost sure
result. The proof of the next result is omitted since it is based on the
relationship \eqref{rel1betweenS} and follows the lines of the proof of
Theorem \ref{stationary} (with obvious modifications). Namely, we prove that
the Stieltjes transform of $\mathbf{B}_{N}$ converges almost surely to the
Stieltjes transform of the same matrix but with the $X_{ij}$'s replaced by
independent real-valued Gaussian random variables with same variance.

\begin{theorem}
\label{M-Pastur}Suppose we have matrices $\mathbf{X}=(X_{ij})_{1\leq i\leq
p,1\leq j\leq n}$ with $(X_{\mathbf{u}})_{\mathbf{u}\in{\mathbb{Z}}^{2}}$ a
strictly stationary real-valued random field. For any $a \in{\mathbb{N}}$, let
$B_{{\mathbf{0}}}^{a}$ be defined in \eqref{defB}, ${\widetilde{\mathcal{F}}%
}^{a}_{\mathbf{0}}=\sigma(X_{\mathbf{u}};\mathbf{u}\in B_{{\mathbf{0}}}^{a})$
and ${\widetilde{\mathcal{F}}}^{\infty}_{\mathbf{0}}=\cap_{a\in\mathbb{N}%
}{\widetilde{\mathcal{F}}}_{\mathbf{0}}^{a}$ (here $\mathbf{0}=(0,0)$). Assume
that $p/n\rightarrow y$ where $y\in(0,\infty)$ and
\[
\mathbb{E}X_{\mathbf{0}}^{2}=1,\text{ }\mathbb{E(}X_{\mathbf{0}}|
{\widetilde{\mathcal{F}}}_{\mathbf{0}}^{1})=0\text{ a.s. and }\mathbb{E(}%
X_{\mathbf{0}}^{2}| {\widetilde{\mathcal{F}}}_{\mathbf{0}}^{\infty})=1\text{
a.s.}
\]
Then,
\[
\mathbf{F}^{\mathbf{XX}^{T}/n}\ \Rightarrow\tilde{G}_{y}\ \text{ a.s.}%
\]
where $\tilde{G}_{y}$ is the standard Marchenko-Pastur distribution function.
\end{theorem}

Note that the above theorem extends the Marchenko-Pastur convergence theorem
from the i.i.d. case to the martingale differences case without additional
moment assumption.

\medskip

3. \textbf{Other convergence results.}

Our results could be also combined with other theorems for Gaussian
structures. If, for instance, the covariance structure is of the form%
\begin{equation}
cov(X_{ij},X_{uv})=a_{i}^{2}a_{j}^{2}I(i=u)I(j=v)+a_{i}^{2}a_{j}%
^{2}I(i=v)I(j=u)\,. \label{cov}%
\end{equation}
with%
\begin{equation}
\max_{j\geq1}|a_{j}|<\infty\,, \label{l2}%
\end{equation}
then%
\[
cov(X_{ij},X_{uv})=V(i,u)V(j,v)+V(i,v)V(j,u)\,,
\]
where $V(i,u)=a_{i}^{2}I(i=u).$ We note that condition (2.1) in Boutet de
Monvet and Khorunzhy (1999) is satisfied and their Theorem 2.2 applies via our
Theorems \ref{replacement} or \ref{almost sure} where we reduced the study to
independent Gaussian variables. This is exactly the function $V(i,u)$ treated
in their Remark (iv) on page 918. The spectral limit can be specified uniquely
by the relations (2.9a) and (2.9b) in Boutet de Monvet and Khorunzhy (1999)
provided the following limit exists%
\begin{equation}
\nu(t)=\lim_{n\rightarrow\infty}\frac{1}{n}\sum_{1\leq j\leq n}I(a_{j}^{2}\leq
t)\,. \label{defniu}%
\end{equation}
More precisely we obtain
\begin{equation}
\label{ConvMonv}\ \mathbf{F}^{{\mathbb{X}}_{n}}\Rightarrow\mathbf{F}\text{
a.s. }%
\end{equation}
where the Stieltjes transform of $\mathbf{F}$ is given by the relation%
\[
S(z)=\int_{0}^{\infty}\frac{d\nu(\lambda)}{-z-\lambda g(z)}\,,
\]
where $g(z)$ is solution of the equation
\[
g(z)=\int_{0}^{\infty}\frac{\lambda d\nu(\lambda)}{-z-\lambda g(z)}\text{
\ \ }z\in{\mathbb{C}}\backslash{\mathbb{R}}\,.
\]
This equation is uniquely solvable in the class of analytic functions $f$
defined on ${\mathbb{C}}\backslash{\mathbb{R}}$ satisfying the conditions
\[
\lim_{x\rightarrow\infty}xf(ix)<\infty,\text{ \ }\operatorname{Im}%
f(z)\operatorname{Im}z>0\text{ for }z\in{\mathbb{C}}\backslash{\mathbb{R}}\,.
\]

Therefore we can formulate the following corollary:

\begin{corollary}
Assume that $(X_{ij})$ are as in Theorem \ref{almost sure} and conditions
(\ref{cov}), (\ref{l2}) and (\ref{defniu}) are satisfied. Then, the
convergence (\ref{ConvMonv}) holds.
\end{corollary}

This result can be applied if $(a_{j}^{2})$ are selected from a stationary and
ergodic sequence of random variables $(A_{k}^{2})$ with distribution function
$\nu(t)$ and such that $|A_{k}|<Y$ $\ $a.s. for some positive random variable
$Y$. In this case, there is a subset $\Omega^{\prime}\subset\Omega,$ with
$\mathbb{P}(\Omega^{\prime})=1$ such that for all $\omega\in\Omega^{\prime}$
\[
\lim_{n\rightarrow\infty}\frac{1}{n}\sum_{1\leq j\leq n}I(A_{j}^{2}\leq
t)(\omega)=\nu(t)\text{ and }|A_{k}(\omega)|<Y(\omega)\,.
\]
Then, for $a_{k}^{2}=A_{k}^{2}(\omega)$, the convergence (\ref{ConvMonv}) holds.

\textbf{ }

\bigskip

\section{Applications}

We mention now three applications of our results to classes of random matrices
with martingale differences entries which could not be treated by the previous
results in the literature. Notice that such results are relevant to
statistical procedures. They give, for instance, theoretical justification to
use the so-called Wachter plot introduced in \cite{Wa}.

\bigskip

\medskip

\textbf{Example 1}. We consider a non linear ARCH($\infty$) random field
$(X_{ij})_{(i,j)\in\mathbb{Z}^{2}}$ given by
\begin{equation}
X_{ij}=\xi_{ij}(c+\sum_{(k,\ell)>_{\text{lex}}(0,0)}g_{k\ell}(X_{i-k,j-\ell
}))\,, \label{ARCHmodel}%
\end{equation}
where $(\xi_{ij})_{(i,j)\in\mathbb{Z}^{2}}$ is a sequence of centered i.i.d.
real-valued random variables such that $\Vert\xi_{\mathbf{0}}\Vert_{2}=1$,
$c>0$ and the $g_{k\ell}$ are functions from ${\mathbb{R}}$ to ${\mathbb{R}}$
such that for any $(x,y)\in{\mathbb{R}}^{2}$,
\[
|g_{k\ell}(x)-g_{k\ell}(y)|\leq\alpha_{k\ell}|x-y|\,.
\]
If $\sum_{(i,j)>_{\text{lex}}(0,0)}\alpha_{ij}<1$ then, by Corollary 2 p. 121
in Doukhan and Truquet (2007), there exists a unique stationary solution of
equation (\ref{ARCHmodel}). This solution is in $\mathbb{L}^{2}$ and can be
written as $X_{ij}=g((\xi_{i-k,j-\ell})_{(k,\ell)\geq_{\text{lex}}(0,0)})$.
Denote $\sigma^{2}=\mathbb{E}(X_{\mathbf{0}}^{2}).$ Based on this stationary
random field we construct the symmetric random matrix ${\mathbb{X}}_{n}.$

For any non-negative integer $a$, consider the sigma algebras ${\mathcal{G}%
}_{ij}^{a}$ and ${\widetilde{\mathcal{F}}}_{ij}^{a}$ defined by
\[
{\mathcal{G}}_{ij}^{a}=\sigma(\xi_{uv}\,:\,(u,v)\in B_{ij}^{a}) \text{ and }
{\widetilde{\mathcal{F}}}_{ij}^{a}=\sigma(X_{uv}\,:\,(u,v)\in B_{ij}^{a}) \,
,
\]
with $B_{ij}^{a}$ defined by (\ref{defB}). Note that ${\mathcal{\tilde{F}}%
}_{ij}^{a}\subseteq{\mathcal{G}}_{ij}^{a}$. Therefore
\[
\mathbb{E}(X_{ij}|{\widetilde{\mathcal{F}}}_{ij}^{1})=(c+\sum_{(k,\ell
)>_{\text{lex}}(0,0)}g_{k\ell}(X_{i-k,j-\ell}))\mathbb{E}(\xi_{ij}%
|{\widetilde{\mathcal{F}}}_{ij}^{1})=0\text{ \ a.s.}%
\]
In addition, since ${\mathcal{G}}_{\mathbf{0}}^{\infty}=\cap_{a\in{\mathbb{N}%
}}{\mathcal{G}}_{\mathbf{0}}^{a}$ is trivial and ${\widetilde{\mathcal{F}}%
}_{\mathbf{0}}^{\infty}:=\cap_{a\in{\mathbb{N}}}{\widetilde{\mathcal{F}}%
}_{\mathbf{0}}^{a}\subseteq{\mathcal{G}}_{\mathbf{0}}^{\infty}$, it follows
that $\mathbb{E}(X_{\mathbf{0}}^{2}|{\widetilde{\mathcal{F}}}_{\mathbf{0}%
}^{\infty})=\sigma^{2}$ a.s. Therefore all the conditions of Theorem
\ref{stationary} and also of Theorem \ref{M-Pastur} are satisfied and
therefore their conclusions hold for ${\mathbb{X}}_{n}/\sigma$.\medskip

\textbf{Example 2.} Consider a real-valued martingale differences sequence
$(D_{i})_{i\geq1}$ adapted to the natural filtrations $\mathcal{F}_{k}%
=\sigma(D_{j},1\leq j\leq k)$, and with finite second moment. Let
$(\gamma_{ij})$ be a matrix of real-valued random variables which are
independent of $(D_{i})_{i\geq0}$ and with finite second moments. Then
construct the symmetric matrix by using the lexicographic order in the
following way:
\begin{align*}
X_{ij}  &  =\gamma_{ij}D_{u(i,j)}\text{ where }u(i,j)=\frac{(i-1)i}{2}+j\text{
for }1\leq j\leq i\leq n;\\
X_{ij}  &  =X_{ji}\text{ for }1\leq i<j\leq n.
\end{align*}
For clarity we sketch below the lower half of this matrix. The rest is
completed by symmetry.
\[
\mathbf{D}^{n}=\left(
\begin{array}
[c]{cccccc}%
\gamma_{11}D_{1} & \  & \dots &  &  & \\
\gamma_{21}D_{2} & \gamma_{22}D_{3} & \  & \dots &  & \\
\gamma_{31}D_{4} & \gamma_{32}D_{5} & \gamma_{33}D_{6} & \dots &  & \\
\gamma_{41}D_{7} & \gamma_{42}D_{8} & \gamma_{43}D_{9} & \gamma_{44}D_{10} &
\dots & \\
\dots &  &  &  &  & \dots\\
\gamma_{n1}D_{1+n(n-1)/2} & \dots &  &  &  & \gamma_{nn}D_{n(n+1)/2}%
\end{array}
\right)
\]
For any non-negative integer $a$, let us introduce the filtrations
\[
\Gamma_{ij}^{a}=\sigma(\gamma_{uv}\,:\,(u,v)\in B_{ij}^{a}) \,,
\]
where $B_{ij}^{a}$ is defined in \eqref{defB}.

The following result is valid.

\begin{corollary}
\label{cor15} Assume that for some positive $\delta$ we have $\sup
_{i}\mathbb{E}|D_{i}|^{2+\delta}<\infty$ and that there is a positive constant
$c$ such that $\sup_{i,j}|\gamma_{ij}|<c$ a.s. Assume also that%
\[
\lim_{a\rightarrow\infty}\limsup_{n\rightarrow\infty}\frac{1}{n}%
{\displaystyle\sum\limits_{a\leq i\leq n}}|\mathbb{E}(D_{i}^{2}|\mathcal{F}%
_{i-a})-\mathbb{E}D_{i}^{2}|=0\text{ a.s.}\,,
\]%
\[
\lim_{a\rightarrow\infty}\limsup_{n\rightarrow\infty}\frac{1}{n^{2}%
}{\displaystyle\sum\limits_{1\leq j\leq i\leq n}}|\mathbb{E}(\gamma_{ij}%
^{2}|\Gamma_{ij}^{a})-\mathbb{E}\gamma_{ij}^{2}|=0\text{ a.s.}%
\]
Then the conclusion of Theorem \ref{almost sure} holds.
\end{corollary}

The proof of this corollary is a consequence of Theorem \ref{almost sure} via
the following remark which uses the proof of Theorem \ref{almost sure} and
Remark \ref{remark4}:

\begin{remark}
\label{remark2} The conclusion of Theorem \ref{almost sure} holds if we
replace condition (\ref{mix a s}) by the following condition: \newline For any
non-negative integer $a$, there is a filtration $\mathcal{K}_{ij}^{a}$
satisfying for any $j \leq i$: $\mathcal{F}_{ij}^{0}\subseteq\mathcal{K}%
_{ij}^{0}$, $\mathcal{K}_{ij}^{a}\subseteq\mathcal{K}_{i j}^{0}$,
$\mathcal{K}_{ij}^{0}\subseteq\mathcal{K}_{i+1, j}^{0}$ and $\mathcal{K}%
_{i-a,j}^{0}\subseteq\mathcal{K}_{ij}^{a}$ for $i \geq a+1$, and such that
\begin{equation}
\lim_{a\rightarrow\infty}\limsup_{n\rightarrow\infty}\frac{1}{n^{2}}%
\sum_{1\leq j\leq i\leq n}|\mathbb{E}(X_{ij}^{2}-\sigma_{ij}^{2}|{\mathcal{K}%
}_{ij}^{a})|=0\ \text{a.s.} \label{mix a.s.1}%
\end{equation}

\end{remark}

\noindent\textbf{Proof of Corollary \ref{cor15}}. To prove the result, we
first introduce the following notations: for any non-negative integer $a$,
let
\[
v(i,j,a)= \frac{(i-1)i}{2}+(j-a){\mathbf{1}}_{j \geq a+1} + (j-1){\mathbf{1}%
}_{1 \leq j \leq a} ,
\]
\[
{\mathcal{G}}^{a}_{ij}=\mathcal{F}_{v(i,j,a)} \ \text{ and } \ \mathcal{K}%
_{ij}^{a}=\Gamma_{ij}^{a}\vee\mathcal{G}^{a}_{ij} \, .
\]
It is easy to see that for any $j \leq i$, the filtration $\mathcal{K}%
_{ij}^{a}$ satisfies the inclusion properties of Remark \ref{remark2}. Now, by
the independence between the sequences $(D_{i})$ and $(\gamma_{ij})$, we have
\[
\mathbb{E}(\gamma_{ij}D_{u(i,j)}|\mathcal{K}_{ij}^{1})=\mathbb{E}(\gamma
_{ij}|\Gamma_{ij}^{1})\mathbb{E(}D_{u(i,j)}|\mathcal{F}_{u(i,j)-1})=0\text{
a.s.}%
\]
According to Theorem \ref{almost sure} and Remark \ref{remark2}, the corollary
will follow if we shall check the condition \eqref{mix a.s.1} for
$\mathcal{K}_{ij}^{a}$ defined above. Simple algebra shows that
\begin{gather*}
|\mathbb{E}[\gamma_{ij}^{2}D_{u(i,j)}^{2}-\mathbb{E}(\gamma_{ij}%
^{2})\mathbb{E}(D_{u(i,j)}^{2})|\mathcal{K}_{ij}^{a}]|\leq\mathbb{E}%
(\gamma_{ij}^{2}|\Gamma_{ij}^{a})|\mathbb{E}(D_{u(i,j)}^{2}|\mathcal{G}%
_{ij}^{a})-\mathbb{E}(D_{u(i,j)}^{2})|\\
+\mathbb{E}(D_{u(i,j)}^{2})|\mathbb{E}(\gamma_{ij}^{2}|\Gamma_{ij}%
^{a})-\mathbb{E}(\gamma_{ij}^{2})|\,.
\end{gather*}
Clearly, under the conditions of Corollary \ref{cor15}, condition
\eqref{mix a.s.1} will hold if we prove that
\begin{equation}
\label{cor15p1}\lim_{a\rightarrow\infty}\limsup_{n\rightarrow\infty}\frac
{1}{n^{2}}{\displaystyle\sum\limits_{1\leq j\leq i\leq n}}|\mathbb{E}%
(D_{u(i,j)}^{2}|\mathcal{G}_{ij}^{a})-\mathbb{E}(D_{u(i,j)}^{2})| =0\text{
a.s.}%
\end{equation}
With this aim, we write
\begin{multline}
\label{cor15p2}\frac{1}{n^{2}}{\displaystyle\sum\limits_{1\leq j\leq i\leq n}%
}|\mathbb{E}(D_{u(i,j)}^{2}|\mathcal{G}_{ij}^{a})-\mathbb{E}(D_{u(i,j)}%
^{2})|\\
\leq\frac{1}{n^{2}}{\displaystyle\sum_{i=a+1}^{n} \sum_{j=a+1}^{i}}%
|\mathbb{E}(D_{u(i,j)}^{2}|\mathcal{F}_{u(i,j)-a})-\mathbb{E}(D_{u(i,j)}%
^{2})|\\
+ \frac{1}{n^{2}}{\displaystyle\sum_{j=1}^{a}\sum_{i=j}^{n} }\mathbb{E}%
(D_{u(i,j)}^{2})+ \frac{1}{n^{2}}{\displaystyle\sum_{j=1}^{a}\sum_{i=j}^{n}
}\mathbb{E}(D_{u(i,j)}^{2}|\mathcal{F}_{u(i,j)-1}) \, .
\end{multline}
By assumption, the first term in the right-hand side is going to zero when we
first let $n$ tend to infinity and after $a$. Clearly the second one is going
to zero as $n$ is going to infinity since we have $\sup_{i} {\mathbb{E}} (
D_{i}^{2}) < \infty$. To handle the third term, we use the following
decomposition:
\[
\frac{1}{n^{2}}{\displaystyle\sum_{j=1}^{a}\sum_{i=j}^{n} }\mathbb{E}%
(D_{u(i,j)}^{2}|\mathcal{F}_{u(i,j)-1}) \leq\frac{a}{n} + \frac{1}%
{n^{2+\delta/2}}{\displaystyle\sum_{j=1}^{a}\sum_{i=j}^{n} }\mathbb{E}%
(|D_{u(i,j)}|^{2+\delta}I(| D_{u(i,j)}|> n^{1/2} ) | \mathcal{F}_{u(i,j)-1})\,
,
\]
where $\delta$ is such that $\sup_{i}\mathbb{E}(|D_{i}|^{2+\delta}) < \infty$.
Since
\[
\sum_{n \geq1} \frac{1}{n^{2+\delta/2}}{\displaystyle\sum_{j=1}^{a}\sum
_{i=j}^{n} }\mathbb{E}(|D_{u(i,j)}|^{2+\delta})< \infty\, ,
\]
we conclude easily that the last term in the right-hand side of
\eqref{cor15p2} converges to zero as $n$ tend to infinity. This ends the proof
of condition \eqref{mix a.s.1} and therefore of the corollary. $\lozenge$

\bigskip

We list below another corollary which follows from our Theorem
\ref{replacement} and whose proof is straightforward.

\begin{corollary}
Assume that $(\gamma_{ij})$ is a sequence of constants satisfying
$\sup_{(i,j)}|\gamma_{ij}|<\infty$ and assume
\begin{equation}
\label{c0cor16}\frac{1}{n}{\sum\limits_{1\leq i\leq n}}\mathbb{E}(D_{i}%
^{2})<\infty\,,
\end{equation}
\begin{equation}
\label{c1cor16}\lim_{a}\limsup_{n}\frac{1}{n}{\sum\limits_{a\leq i\leq n}%
}\Vert\mathbb{E}(D_{i}^{2}|\mathcal{F}_{i-a})-\mathbb{E}D_{i}^{2}\Vert
_{1}=0\,,
\end{equation}
and for any $\varepsilon>0$,
\begin{equation}
\label{c2cor16}\frac{1}{n^{2}}\sum_{1\leq i\leq n^{2}}\mathbb{E}(D_{i}%
^{2}I(|D_{i}|>\varepsilon n^{1/2}))\rightarrow0\,.
\end{equation}
Then the conclusion of Theorem \ref{replacement} holds.
\end{corollary}

In particular, if the sequence $(\gamma_{ij})$ is constant, the only relevant
conditions in these two last corollaries are imposed on the differences of
martingale. Notice also that if $(D_{i}, i \in{\mathbb{Z}})$ is a strictly
stationary sequence of martingale differences in ${\mathbb{L}}^{2}$, the
conditions \eqref{c0cor16} and \eqref{c2cor16} are obviously satisfied and
\eqref{c1cor16} becomes ${\mathbb{E}}(D_{0}^{2} | {\mathcal{F}}_{-\infty})=
{\mathbb{E}}(D_{0}^{2}) $ in ${\mathbb{L}}^{1}$, where ${\mathcal{F}}%
_{-\infty} = \cap_{i \in{\mathbb{Z}}} \sigma( D_{k},k\leq i)$. This last
condition is equivalent to ${\mathbb{E}}(D_{0}^{2} | {\mathcal{F}}_{-\infty})=
{\mathbb{E}}(D_{0}^{2})$ a.s. and it holds if the sequence is ergodic or
strong mixing.

\medskip

\textbf{Example 3.} Consider $p$ independent copies $(D^{(i)}_{j})_{j
\in{\mathbb{Z}}}$, $i=1, \dots, p$ of a real-valued martingale differences
sequence $(D_{i})_{i \in{\mathbb{Z}}}$ with respect to the natural filtration
${\mathcal{F}}_{j}=\sigma(D_{k}, k \leq j)$, such that ${\mathbb{E}}(D_{i}%
^{2})=1$ for any $i \in{\mathbb{Z}}$. Let $D_{ij}=D_{j}^{(i)}$ and
$\mathbf{X}=\mathbf{X}_{np} =(D_{ij})_{1\leq i\leq p,1\leq j\leq n}$. Applying
our Theorem \ref{Th-Ma-Pa}, the following corollary holds for the sample
covariance matrix:

\begin{corollary}
\label{directcor} Assume that conditions \eqref{c1cor16} and \eqref{c2cor16}
hold, and that $p/n\rightarrow y \in(0,\infty)$. Then $\mathbf{F}%
^{\mathbf{XX}^{T}/n}\ \Rightarrow\tilde{G}_{y}$ a.s., where $\tilde{G}_{y}$ is
the standard Marchenko-Pastur distribution function.
\end{corollary}

\noindent\textbf{Proof of Corollary \ref{directcor}.} By using the fact that
for any $i\in\{2,\dots,p\}$, $\sigma((D_{j}^{(i)})_{j\in{\mathbb{Z}}})$ is
independent of $\sigma(D_{j}^{(k)})_{j\in{\mathbb{Z}}},1\leq k\leq i-1)$, we
can easily verify that all the conditions of Theorem \ref{Th-Ma-Pa} are
satisfied under the assumptions of Corollary \ref{directcor}. Therefore,
setting ${\mathbf{A}}_{n}=n^{-1}\mathbf{X}\mathbf{X}^{T}$ we obtain
$\mathbf{F}^{{\mathbf{A}}_{n}}\ \Rightarrow\tilde{G}_{y}$ in probability, or
equivalently, for any $z\in{\mathbb{C}}^{+}$, $S^{{\mathbb{A}}_{n}%
}(z)\rightarrow S_{y}(z)$ in probability, where $S_{y}(z)$ is the Stieltjes
transform of $\tilde{G}_{y}$. Furthermore, since both Stieljes transforms are
bounded, the convergence in probability implies ${\mathbb{E}}(S^{{\mathbb{A}%
}_{n}}(z))\rightarrow S_{y}(z).$ Now, since the rows of $\mathbf{X}$ are
independent, for any $z\in{\mathbb{C}}^{+}$ we obtain $S^{{\mathbb{A}}_{n}%
}(z)-{\mathbb{E}}(S^{{\mathbb{A}}_{n}}(z))\rightarrow0$ a.s. (see, for
instance, Lemma 4.1 in \cite{A}). So, overall, under the conditions of
Corollary \ref{directcor}, we get that $S^{{\mathbb{A}}_{n}}(z)$ converges
almost surely to $S_{y}(z)$ that is equivalent to $\mathbf{F}^{{\mathbf{A}%
}_{n}}\ \Rightarrow\tilde{G}_{y}$ a.s. $\lozenge$

\section{Proofs\textbf{ }}

\subsection{Proof of \textbf{Theorem \ref{replacement}}}

\qquad We start this section with some notations. For a function $f$ of one
variable $x$ we denote by $d^{i}f=d^{i}f/dx^{i},$ the derivative of order $i$
with respect to $x.$ For a multivariate function we use the notations
$\partial_{k}^{i}f$ $=\partial^{i}f/\partial^{i}x_{k}$ for the partial
derivative of order $i$ with respect to the variable $x_{k}.$ Also
$\partial_{jk}^{2}f=\partial^{2}f/\partial x_{j}\partial x_{k}$ means the
derivatives with respect to $x_{j}$ of the derivative with respect to $x_{k}$,
and so on.

\medskip

Let $k_{n}=n(n+1)/2$ and $\mathbf{x}_{n}=(x_{ij})_{1\leq j\leq i\leq n}$ be a
vector of ${\mathbb{R}}^{k_{n}}$. Let $\mathbf{A}_{n}(\mathbf{x}_{n})$ be the
symmetric matrix of order $n$ defined by
\begin{equation}
(\mathbf{A}_{n}(\mathbf{x}_{n}))_{ij}=\left\{
\begin{array}
[c]{ll}%
\frac{1}{\sqrt{n}}x_{ij}\  & 1\leq j\leq i\leq n\\
\frac{1}{\sqrt{n}}x_{ji}\  & 1\leq i<j\leq n\,.
\end{array}
\right.  \label{defofA}%
\end{equation}
It is convenient to introduce a function notation for the Stieltjes transform
defined in (\ref{StTr}). Let $z\in\mathbb{C}^{+}$ and let $s_{n}%
(\mathbf{x}_{n},z)$ be the function defined from $\mathbb{R}^{k_{n}}$ to
$\mathbb{C}$ by
\begin{equation}
s(\mathbf{x}_{n})=s_{n}(\mathbf{x}_{n},z)=S^{\mathbf{A}_{n}(\mathbf{x}_{n}%
)}(z)=\frac{1}{n}\mathrm{Tr}(\mathbf{A}_{n}(\mathbf{x}_{n})-z{\mathbf{{I}}%
_{n}})^{-1}\,, \label{defff}%
\end{equation}
where $\mathbf{{I}}_{n}$ is the identity matrix of order $n$ and for
simplicity here and in the sequel we deleted the variable $z$ and the index
$n$ from the notation of $s_{n}(\mathbf{x}_{n},z)$. So we write $s(\mathbf{x}%
_{n})$ instead of $s_{n}(\mathbf{x}_{n},z)$ when no confusion is possible. The
partial derivatives of the function $s(\mathbf{x}_{n})$ have been estimated in
Chatterjee (2006). There are three positive constants $c_{1},c_{2}$ and
$c_{3}$ depending on $\operatorname{Im}z$ such that
\begin{equation}
|\partial_{u}s(\mathbf{x}_{n})|\leq\frac{c_{1}}{n^{3/2}}\text{; }%
|\partial_{uv}^{2}s(\mathbf{x}_{n})|\leq\frac{c_{2}}{n^{2}}\text{ and
}|\partial_{uvw}^{3}s(\mathbf{x}_{n})|\leq\frac{c_{3}}{n^{5/2}}\text{ for all
}u,v,w\text{.} \label{boundsd}%
\end{equation}

The proof is based on Proposition \ref{approximation} given in Section
\ref{TehRes}. We shall order the indexes of the variables $(X_{ij})_{1\leq
j\leq i\leq n}$ by using the lexicographic order. These indexes are denoted by
$u_{1}<_{\mathrm{\ell ex}}u_{2}<_{\mathrm{\ell ex}}\cdots<_{\mathrm{\ell ex}%
}u_{k_{n}}.$ Here is the enumeration for the indexes in the lower half part of
the matrix
\[
\left(
\begin{array}
[c]{ccccccc}%
u_{1} & \dots &  &  &  &  & \\
u_{2} & u_{3} & \dots &  &  & \  & \\
u_{4} & u_{5} & u_{6} & \dots &  &  & \\
u_{7} & u_{8} & u_{9} & u_{10}\  & \dots &  & \\
u_{11} & u_{12} & u_{13} & u_{14} & u_{15} &  & \\
\dots &  &  &  &  &  & \dots\\
u_{1+n(n-1)/2} & \dots &  &  &  &  & u_{k_{n}}%
\end{array}
\right)  \,.
\]
With the above notations, we have that $S^{{\mathbb{X}}_{n}}-S^{{\mathbb{Y}%
}_{n}}=s(\mathbf{X}_{n})-s(\mathbf{Y}_{n})$ where $\mathbf{X}_{n}=(X_{u_{\ell
}})_{1\leq\ell\leq k_{n}}$ and $\mathbf{Y}_{n}=(Y_{u_{\ell}})_{1\leq\ell\leq
k_{n}}$. To prove the theorem, we shall show in what follows that
\begin{equation}
\lim_{n\rightarrow\infty}{\mathbb{E}}|s(\mathbf{X}_{n})-s(\mathbf{Y}%
_{n})|=0\,. \label{convdelta1}%
\end{equation}
We start the proof by truncating the random variables. Let $\varepsilon>0$.
For any integer $\ell\in\lbrack1,k_{n}]$, we then define
\begin{equation}
T_{u_{\ell}}=T_{n,u_{\ell}}=X_{u_{\ell}}I(|X_{u_{\ell}}|\leq\varepsilon
\sqrt{n})\,. \label{defTul}%
\end{equation}
As in our previous notation, when no confusion is possible, to ease the
notation, we shall use the notation of $T_{u_{\ell}}$ instead of
$T_{n,u_{\ell}},$ but we shall keep always in mind the dependence of $n$.
Since $\mathbb{E(}X_{u_{\ell}}-T_{u_{\ell}})^{2}\leq\mathbb{E}(X_{u_{\ell}%
}^{2}I(|X_{u_{\ell}}|>\varepsilon\sqrt{n}))$, by (\ref{L}), it follows that
\begin{equation}
\frac{1}{n^{2}}{\displaystyle\sum\limits_{1\leq\ell\leq k_{n}}}\mathbb{E(}%
X_{u_{\ell}}-T_{u_{\ell}})^{2}\rightarrow0\,\text{ as $n\rightarrow\infty\,.$}
\label{negl1}%
\end{equation}
Denoting $\mathbf{T}_{n}=(T_{n,u_{\ell}})_{1\leq\ell\leq k_{n}}$, and using
Lemma \ref{lmagotze} from Section \ref{TehRes}, we have
\[
{\mathbb{E}}|s(\mathbf{X}_{n})-s(\mathbf{T}_{n})|^{2}\ll\frac{1}{n^{2}%
}{\displaystyle\sum\limits_{1\leq\ell\leq k_{n}}}\mathbb{E(}X_{u_{\ell}%
}-T_{u_{\ell}})^{2}\,.
\]
Taking into account (\ref{negl1}), it follows that
\begin{equation}
{\mathbb{E}}|s(\mathbf{X}_{n})-s(\mathbf{T}_{n})|^{2}\rightarrow0\,\text{as
}n\rightarrow\infty. \label{2}%
\end{equation}
Let us consider now a vector $\mathbf{Z}_{n}=(Z_{n,u_{\ell}})_{1\leq\ell\leq
k_{n}}$ of independent centered real-valued Gaussian random variables,
independent of $\mathbf{X}_{n}$ and such that, for all $\ell$, we have
$\mathbb{E}Z_{n,u_{\ell}}^{2}=\mathbb{E}(X_{u_{\ell}}^{2}I(|X_{u_{\ell}}%
|\leq\varepsilon\sqrt{n}))$. We denote for short $Z_{n,u_{\ell}}=Z_{u_{\ell}%
}.$ Let ${\mathbb{Z}}_{n}$ be the matrix constructed as in (\ref{defX}). By
Lemma \ref{inter} and (\ref{boundsd}), we get that
\[
|\mathbb{E}(s(\mathbf{Y}_{n}))-\mathbb{E}(s(\mathbf{Z}_{n}))| \ll\frac
{1}{n^{2}}{\displaystyle\sum\limits_{1\leq\ell\leq k_{n}}}|\mathbb{E}%
X_{u_{\ell}}^{2}-\mathbb{E}Z_{u_{\ell}}^{2}|=\frac{1}{n^{2}}{\displaystyle\sum
\limits_{1\leq\ell\leq k_{n}}}\mathbb{E}(X_{u_{\ell}}^{2}I(|X_{u_{\ell}%
}|>\varepsilon\sqrt{n}))\,.
\]
Hence, using (\ref{L}), it follows that
\begin{equation}
\lim_{n\rightarrow\infty}|\mathbb{E}(s(\mathbf{Y}_{n}))-\mathbb{E}%
(s(\mathbf{Z}_{n}))|=0\,. \label{negl3}%
\end{equation}
Since $\mathbf{Y}_{n}$ and $\mathbf{Z}_{n}$ have independent components, it is
well-known (see for instance the proof on page 34 in Bai-Silverstein, 2010)
that $s(\mathbf{Y}_{n})-\mathbb{E}s(\mathbf{Y}_{n})\rightarrow0$ a.s. and also
$s(\mathbf{Z}_{n})-\mathbb{E}s(\mathbf{Z}_{n})\rightarrow0$ a.s. By combining
these last two almost sure convergence results with \eqref{negl3}, we get that
$s(\mathbf{Y}_{n})-s({\mathbf{Z}}_{n})\rightarrow0$ a.s. Since the Stieltjes
transforms are bounded, we also derive that
\begin{equation}
\mathbb{E}|s(\mathbf{Y}_{n})-s({\mathbf{Z}}_{n})|\rightarrow0\,\text{ as
$n\rightarrow\infty$ .} \label{4}%
\end{equation}
Therefore by (\ref{2}) and (\ref{4}), we note that the convergence
\eqref{convdelta1}, and then the conclusion of the theorem will follow if we
prove that
\begin{equation}
\lim_{\varepsilon\rightarrow0}\limsup_{n\rightarrow\infty}\mathbb{E}%
|s({\mathbf{T}}_{n})-s({\mathbf{Z}}_{n})|=0\,. \label{neglaimth1}%
\end{equation}
With this aim, we shall apply the approximation in Proposition
\ref{approximation}. Let $a$ be a fixed but arbitrary positive integer. For
$u_{\ell}=(i,j),$ $1\leq j\leq i\leq n$, let $B_{u_{\ell}}(a)$ be the set
\begin{equation}
B_{u_{\ell}}(a)=\{(u,v)\in{\mathbb{N}}^{2}\,:\,v\leq u\,\text{ and }\,(u,v)\in
B_{ij}^{1}\setminus B_{ij}^{a}\} \label{defB(u,a)}%
\end{equation}
with $B_{ij}^{1}$ and $B_{ij}^{a}$ defined by (\ref{defB}). For example, when
$a=2$, the indexes $(u,v)$ that belong to the set $B_{u_{\ell}}(2)$ are
described by the points in the next matrix (below $u_{\ell}=(i,j)$ with $i\geq
j+2$).
\[
\left(
\begin{array}
[c]{cccccccccccccc}
&  &  &  &  &  &  &  &  &  &  &  &  & \\
&  &  &  &  &  &  &  &  &  &  &  &  & \\
&  &  &  &  &  &  &  &  &  &  &  &  & \\
&  & . & . & . &  &  &  &  &  &  &  &  & \\
&  & . & u_{\ell} &  &  &  &  &  &  &  &  &  & \\
&  &  &  &  &  &  &  &  &  &  &  &  &
\end{array}
\right)
\]
Denote
\[
\mathbf{C}_{u_{\ell}}=\mathbf{C}_{n,u_{\ell}}=(T_{u_{1}},\dots,T_{u_{\ell-1}%
},0,Z_{u_{\ell+1}},\dots,Z_{u_{k_{n}}})
\]
and
\[
\mathbf{U}_{u_{\ell-1}}=\mathbf{U}_{n,u_{\ell-1}}=\mathbf{(}U_{u_{1}}%
,\dots,U_{u_{\ell-1}})\,,
\]
where $U_{u_{i}}=0$ if $u_{i}\in B_{u_{\ell}}(a)$ and $U_{u_{i}}=T_{u_{i}}$ if
$u_{i}\in B_{u_{\ell}}^{a}$. Set also
\[
\widetilde{\mathbf{C}}_{u_{\ell}}^{a}=\widetilde{\mathbf{C}}_{n,u_{\ell}}%
^{a}=(\mathbf{U}_{u_{\ell-1}},0,Z_{u_{\ell+1}},\dots,Z_{u_{k_{n}}})\,.
\]

\noindent By Proposition \ref{approximation} and (\ref{boundsd}), we get that
\[
s(\mathbf{T}_{n})-s(\mathbf{Z}_{n}):=\mathbb{\ }R_{1,n}+\mathbb{\ }%
R_{2,n}(a)+\mathbb{\ }R_{3,n}(a)\,,
\]
where
\begin{equation}
R_{1,n}={\displaystyle\sum\limits_{1\leq\ell\leq k_{n}}}\mathbb{(}T_{u_{\ell}%
}-Z_{u_{\ell}})\partial_{u_{\ell}}s(\mathbf{C}_{u_{\ell}})\,, \label{defR1n}%
\end{equation}%
\begin{equation}
R_{2,n}(a)=\frac{1}{2}{\displaystyle\sum\limits_{1\leq\ell\leq k_{n}}%
}\mathbb{(}T_{u_{\ell}}^{2}-Z_{u_{\ell}}^{2})\partial_{u_{\ell}}%
^{2}s(\widetilde{\mathbf{C}}_{u_{\ell}}^{a})\,, \label{defR2n}%
\end{equation}
and
\begin{equation}
\label{defR3n}|R_{3,n}(a)|\leq c_{3}\frac{1}{n^{5/2}}{\displaystyle\sum
\limits_{1\leq\ell\leq k_{n}}}\mathbb{(}T_{u_{\ell}}^{2}+Z_{u_{\ell}}%
^{2}){\displaystyle\sum\limits_{u_{k}\in B_{u_{\ell}}(a)}}|T_{u_{k}}|
+c_{3}\frac{1}{n^{5/2}}{\displaystyle\sum\limits_{1\leq\ell\leq k_{n}}%
}(|T_{u_{\ell}}|^{3}+|Z_{u_{\ell}}|^{3})\,.
\end{equation}
We first handle the term $R_{1,n}$ and we write
\begin{equation}
|R_{1,n}|\leq\Big |{\displaystyle\sum\limits_{1\leq\ell\leq k_{n}}}T_{u_{\ell
}}\partial_{u_{\ell}}s(\mathbf{C}_{u_{\ell}})\Big |+\Big |{\displaystyle\sum
\limits_{1\leq\ell\leq k_{n}}}Z_{u_{\ell}}\partial_{u_{\ell}}s(\mathbf{C}%
_{u_{\ell}})\Big |\,. \label{decR1n*}%
\end{equation}
Since the r.v.'s $Z_{u_{\ell}}\partial_{u_{\ell}}s(\mathbf{C}_{u_{\ell}})$,
$1\leq\ell\leq k_{n}$, are orthogonal, by using \eqref{boundsd}, we get
\[
\mathbb{E}\Big |{\displaystyle\sum\limits_{1\leq\ell\leq k_{n}}}Z_{u_{\ell}%
}\partial_{u_{\ell}}s(\mathbf{C}_{u_{\ell}})\Big |^{2}\ll\frac{1}{n^{3}%
}{\displaystyle\sum\limits_{1\leq\ell\leq k_{n}}}\mathbb{E}(Z_{u_{\ell}}%
^{2})\leq\frac{1}{n^{3}}{\displaystyle\sum\limits_{1\leq\ell\leq k_{n}}%
}\mathbb{E}(X_{u_{\ell}}^{2})\,.
\]
Then, by (\ref{bound}), it follows that
\begin{equation}
\mathbb{E}\Big |{\displaystyle\sum\limits_{1\leq\ell\leq k_{n}}}Z_{n,u_{\ell}%
}\partial_{u_{\ell}}s(\mathbf{C}_{n,u_{\ell}})\Big |^{2}\rightarrow0\,,\text{
as $n\rightarrow\infty$}. \label{decR1nbis}%
\end{equation}
To analyze the first term in the right-hand side of \eqref{decR1n*} we use the
following decomposition:
\[
T_{u_{\ell}}\partial_{u_{\ell}}s(\mathbf{C}_{u_{\ell}})=D_{u_{\ell}}%
\partial_{u_{\ell}}s(\mathbf{C}_{u_{\ell}})+{\mathbb{E}}(T_{u_{\ell}%
}|{{\mathcal{F}}_{u_{\ell}}^{1}})\partial_{u_{\ell}}s(\mathbf{C}_{u_{\ell}%
})\,,
\]
where
\[
D_{u_{\ell}}=D_{n,u_{\ell}}=T_{u_{\ell}}-{\mathbb{E}}(T_{u_{\ell}%
}|{{\mathcal{F}}_{u_{\ell}}^{1}})\,.
\]
By using the fact that $\mathbb{E}(X_{u_{\ell}}|\mathcal{F}_{u_{\ell}}^{1})=0$
a.s. and \eqref{boundsd}, we get
\[
\Big |{\displaystyle\sum\limits_{1\leq\ell\leq k_{n}}}{\mathbb{E}}(T_{u_{\ell
}}|{{\mathcal{F}}_{u_{\ell}}^{1}})\partial_{u_{\ell}}s(\mathbf{C}_{u_{\ell}%
})\Big |\ll\frac{1}{n^{3/2}}{\displaystyle\sum\limits_{1\leq\ell\leq k_{n}}%
}{\mathbb{E}}(|X_{u_{\ell}}|I(|X_{u_{\ell}}|>\varepsilon\sqrt{n}%
)|{{\mathcal{F}}_{u_{\ell}}^{1}})\,.
\]
Therefore, by condition \eqref{L},
\[
{\mathbb{E}}\Big |{\displaystyle\sum\limits_{1\leq\ell\leq k_{n}}}{\mathbb{E}%
}(T_{u_{\ell}}|{{\mathcal{F}}_{u_{\ell}}^{1}})\partial_{u_{\ell}}%
s(\mathbf{C}_{u_{\ell}})\Big |
\ll\frac{1}{\varepsilon n^{2}}{\displaystyle\sum\limits_{1\leq\ell\leq k_{n}}%
}{\mathbb{E}}(X_{u_{\ell}}^{2}I(|X_{u_{\ell}}|>\varepsilon\sqrt{n}%
))\rightarrow0\,,\text{ as $n\rightarrow\infty$}.
\]
On the other hand, since the r.v.'s $D_{u_{\ell}}\partial_{u_{\ell}%
}s(\mathbf{C}_{u_{\ell}})$, $1\leq\ell\leq k_{n}$, are orthogonal, by using
\eqref{boundsd}, we get
\[
\mathbb{E}\Big |{\displaystyle\sum\limits_{1\leq\ell\leq k_{n}}}D_{u_{\ell}%
}\partial_{u_{\ell}}s(\mathbf{C}_{u_{\ell}})\Big |^{2}\ll\frac{1}{n^{3}%
}{\displaystyle\sum\limits_{1\leq\ell\leq k_{n}}}\mathbb{E}(D_{u_{\ell}}%
^{2})\,.
\]
But, by the properties of the conditional expectation, $\mathbb{E}(D_{u_{\ell
}}^{2})\leq\mathbb{E}(T_{u_{\ell}}^{2})\leq\mathbb{E}(X_{u_{\ell}}^{2})$.
Hence, by using (\ref{bound}), it follows that
\[
\mathbb{E}\Big |{\displaystyle\sum\limits_{1\leq\ell\leq k_{n}}}D_{n,u_{\ell}%
}\partial_{u_{\ell}}s(\mathbf{C}_{n,u_{\ell}})\Big |^{2}\rightarrow0\,,\text{
as $n\rightarrow\infty$}.
\]
So, overall,
\[
\mathbb{E}\Big |{\displaystyle\sum\limits_{1\leq\ell\leq k_{n}}}T_{n,u_{\ell}%
}\partial_{u_{\ell}}s(\mathbf{C}_{n,u_{\ell}})\Big |\rightarrow0\text{ \ as
$n\rightarrow\infty$},
\]
which combined with \eqref{decR1nbis} proves that
\[
\mathbb{E}|R_{1,n}|\rightarrow0\text{ as $\ n\rightarrow\infty$}.
\]

We estimate now the term $\mathbb{E}|R_{3,n}(a)|$. We first note that the
cardinality of $B_{u_{\ell}}(a)$ is smaller than $b=2a(a-1)\leq2a^{2}$.
Therefore, by the level of truncation, we derive
\[
{\displaystyle\sum\limits_{u_{i}\in B_{u_{\ell}}(a)}}|T_{u_{i}}|\leq
2a^{2}\varepsilon\sqrt{n}\,.
\]
Moreover $\mathbb{E}(T_{u_{\ell}}^{2}+Z_{u_{\ell}}^{2})\leq2\sigma_{u_{\ell}%
}^{2}$ and $\mathbb{E}|T_{u_{\ell}}|^{3}\leq\varepsilon n^{1/2}\sigma
_{u_{\ell}}^{2}$. On another hand, since $Z_{u_{\ell}}$ is a Gaussian r.v., it
follows that
\[
\ \mathbb{E}|Z_{u_{\ell}}|^{3}\leq2(\mathbb{E}Z_{u_{\ell}}^{2})^{3/2}%
=2(\mathbb{E}(X_{u_{\ell}}^{2}I(|X_{u_{\ell}}|\leq\varepsilon n^{1/2}%
)))^{3/2}\leq2\varepsilon n^{1/2}\sigma_{u_{\ell}}^{2}\,.
\]
Therefore, the above considerations show that
\[
\mathbb{E}|R_{3,n}(a)|\ll\frac{a^{2}\varepsilon}{n^{2}}{\displaystyle\sum
\limits_{1\leq\ell\leq k_{n}}}\sigma_{u_{\ell}}^{2}\,.
\]
Whence by (\ref{bound}), for any positive integer $a$,
\[
\lim_{\varepsilon\rightarrow0}\limsup_{n\rightarrow\infty}\mathbb{E}%
|R_{3,n}(a)|=0\,.
\]
It remains to analyze $\mathbb{E}|R_{2,n}(a)|$. We shall use the following
decomposition:
\begin{multline}
2R_{2,n}(a)={\displaystyle\sum\limits_{1\leq\ell\leq k_{n}}}\mathbb{(}%
T_{u_{\ell}}^{2}-\mathbb{E(}T_{u_{\ell}}^{2}|\mathcal{F}_{u_{\ell}}%
^{a}))\partial_{u_{\ell}}^{2}s(\widetilde{\mathbf{C}}_{u_{\ell}}%
^{a})\label{decI1I2}\\
+{\displaystyle\sum\limits_{1\leq\ell\leq k_{n}}}\mathbb{(E(}T_{u_{\ell}}%
^{2}|\mathcal{F}_{u_{\ell}}^{a})-Z_{u_{\ell}}^{2})\partial_{u_{\ell}}%
^{2}s(\widetilde{\mathbf{C}}_{u_{\ell}}^{a}):=I_{n}(a)+I\!I_{n}(a).
\end{multline}
The analysis of $I_{n}(a)$ is tedious and is based on a blocking technique
which introduces martingale structure. The estimate is done in Lemma
\ref{moment computation} of Section \ref{TehRes}, which we shall use with
$p=2$, $K=\varepsilon n^{1/2}$, $A_{u_{\ell}}=\partial_{u_{\ell}}%
^{2}s(\widetilde{\mathbf{C}}_{u_{\ell}}^{a})$ and ${\mathcal{G}}%
=\sigma({\mathbf{Z}_{n}})$. Note that by \eqref{boundsd}, $\max_{1\leq\ell\leq
k_{n}}|A_{u_{\ell}}|\leq c_{2}n^{-2}$. It follows that, for any positive
integers $n$ and $a$,
\begin{equation}
\mathbb{E}|I_{n}(a)|\ll\varepsilon\sqrt{a}\,. \label{convL1In}%
\end{equation}
To handle the second term $I\!I_{n}(a)$ in \eqref{decI1I2}, we first apply the
triangle inequality and use \eqref{boundsd} to get
\begin{equation}
\label{termIIn}|I\!I_{n}(a)|\leq\frac{c_{2}}{n^{2}}{\displaystyle\sum
\limits_{1\leq\ell\leq k_{n}}}|\mathbb{E(}T_{u_{\ell}}^{2}|\mathcal{F}%
_{u_{\ell}}^{a})-\mathbb{E}Z_{u_{\ell}}^{2}| +|{\displaystyle\sum
\limits_{1\leq\ell\leq k_{n}}}\mathbb{(}Z_{u_{\ell}}^{2}-\mathbb{E}Z_{u_{\ell
}}^{2})\partial_{u_{\ell}}^{2}s(\widetilde{\mathbf{C}}_{u_{\ell}}^{a})|\,.
\end{equation}
Note that $\mathbb{E}Z_{u_{\ell}}^{2}=\mathbb{E}T_{u_{\ell}}^{2}$. Therefore
\begin{gather*}
\frac{1}{n^{2}}{\displaystyle\sum\limits_{1\leq\ell\leq k_{n}}}\Vert
\mathbb{E}(T_{u_{\ell}}^{2}|\mathcal{F}_{u_{\ell}}^{a})-\mathbb{E}(Z_{u_{\ell
}}^{2})\Vert_{1}\\
\leq\frac{1}{n^{2}}{\displaystyle\sum\limits_{1\leq\ell\leq k_{n}}}%
\Vert\mathbb{E}(X_{u_{\ell}}^{2}|\mathcal{F}_{u_{\ell}}^{a})-\mathbb{E}%
(X_{u_{\ell}}^{2})\Vert_{1}+\frac{2}{n^{2}}{\displaystyle\sum\limits_{1\leq
\ell\leq k_{n}}}\mathbb{E}(X_{u_{\ell}}^{2}I(|X_{u_{\ell}}|>\varepsilon
\sqrt{n}))\,.
\end{gather*}
Hence, taking into account conditions (\ref{mix}) and (\ref{L}), it follows
that
\begin{equation}
\lim_{a\rightarrow\infty}\limsup_{\varepsilon\rightarrow0}\limsup
_{n\rightarrow\infty}\frac{1}{n^{2}}{\displaystyle\sum\limits_{1\leq\ell\leq
k_{n}}}\Vert\mathbb{E}(T_{n,u_{\ell}}^{2}|\mathcal{F}_{u_{\ell}}%
^{a})-\mathbb{E}Z_{n,u_{\ell}}^{2}\Vert_{1}=0\,. \label{convL1Ta}%
\end{equation}
We handle now the last term in the right-hand side of \eqref{termIIn}. Set
\[
d_{u_{\ell}}^{\prime}=d_{n,u_{\ell}}^{\prime}=\mathbb{(}Z_{u_{\ell}}%
^{2}-\mathbb{E}Z_{u_{\ell}}^{2})\partial_{u_{\ell}}^{2}s(\widetilde
{\mathbf{C}}_{u_{\ell}}^{a})\,,
\]
and observe that the r.v.'s $(d_{u_{\ell}}^{\prime})_{\ell\geq1}$ are
orthogonal. Therefore, by (\ref{boundsd}),
\[
\mathbb{E}|{\displaystyle\sum\limits_{1\leq\ell\leq k_{n}}}d_{u_{\ell}%
}^{\prime}|^{2}\leq{\displaystyle\sum\limits_{1\leq\ell\leq k_{n}}}%
\mathbb{E}|d_{u_{\ell}}^{\prime}|^{2}\ll\frac{1}{n^{4}}{\displaystyle\sum
\limits_{1\leq\ell\leq k_{n}}}\mathbb{E}|Z_{u_{\ell}}|^{4}\,.
\]
But by the definition of $Z_{n,u_{\ell}},$ we have%
\[
\mathbb{E}|Z_{u_{\ell}}|^{4}\leq3(\mathbb{E}Z_{u_{\ell}}^{2})^{2}%
=3(\mathbb{E}(X_{u_{\ell}}^{2}I(|X_{u_{\ell}}|\leq\varepsilon n^{1/2}%
)))^{2}\leq3\varepsilon^{2}n\sigma_{u_{\ell}}^{2}\,.
\]
So, by (\ref{bound}),
\begin{equation}
\mathbb{E}|{\displaystyle\sum\limits_{1\leq\ell\leq k_{n}}}d_{u_{\ell}%
}^{\prime}|^{2}\ll\frac{\varepsilon^{2}n}{n^{4}}{\displaystyle\sum
\limits_{1\leq\ell\leq k_{n}}}\sigma_{u_{\ell}}^{2}\rightarrow0\,\ \text{as
}n\rightarrow\infty\,. \label{convL1Ta*}%
\end{equation}
Therefore, from \eqref{convL1Ta} and \eqref{convL1Ta*}, it follows that
\[
\lim_{a\rightarrow\infty}\limsup_{\varepsilon\rightarrow0}\limsup
_{n\rightarrow\infty}\mathbb{E}|I\!I_{n}(a)|=0\,.
\]
Hence, letting $\varepsilon$ tend to zero in \eqref{convL1In}, we get that
\[
\lim_{a\rightarrow\infty}\limsup_{\varepsilon\rightarrow0}\limsup
_{n\rightarrow\infty}\mathbb{E}|R_{2,n}(a)|=0\,.
\]
This ends the proof of the theorem. $\lozenge$

\subsection{Proof of Theorem \ref{almost sure}}

We shall use the same notations as those introduced in the proof of Theorem
\ref{replacement}, and we also start with a truncation argument. For any
integer $\ell$ belonging to $[1,k_{n}]$, let $T_{n,u_{\ell}}$ be defined as in
\eqref{defTul} but with $\varepsilon=1$. Therefore, all along the proof, we
set
\begin{equation}
T_{u_{\ell}}:=T_{n,u_{\ell}}=X_{u_{\ell}}I(|X_{u_{\ell}}|\leq\sqrt{n})\,,
\label{defTul2}%
\end{equation}
$\mathbf{T}_{n}=(T_{n,u_{\ell}})_{1\leq\ell\leq k_{n}}$ and $\mathbf{X}%
_{n}=(X_{u_{\ell}})_{1\leq\ell\leq k_{n}}$. In the rest of the proof, we shall
write $T_{u_{\ell}}$ instead of $T_{n,u_{\ell}}$ when no confusion is
possible. We start by proving that
\begin{equation}
|s({\mathbf{X}}_{n})-s({\mathbf{T}}_{n})|\rightarrow0\,\text{ a.s. as
}n\rightarrow\infty\,\text{.} \label{aim1prth2}%
\end{equation}
By Lemma \ref{lmagotze}, if $z=u+iv$ with $v>0$,
\begin{align}
|s({\mathbf{X}}_{n})-s({\mathbf{T}}_{n})|^{2}  &  \leq\frac{1}{n^{2}v^{4}%
}{\sum\limits_{1\leq\ell\leq k_{n}}}(X_{u_{\ell}}-{T}_{u_{\ell}}%
)^{2}\nonumber\label{appligotze1}\\
&  \leq\frac{1}{n^{2}v^{4}}{\sum\limits_{1\leq\ell\leq k_{n}}}X_{u_{\ell}}%
^{2}I(|X_{u_{\ell}}|>n^{1/2}):=v^{-4}U_{n}\,.
\end{align}
Hence, by the Borel-Cantelli lemma, in order to prove \eqref{aim1prth2}, it is
enough to prove that, for any $\varepsilon>0$,
\[
\sum_{r\geq0}\mathbb{P}\big (\max_{2^{r}\leq j<2^{r+1}}U_{j}>\varepsilon
\big )<\infty\,.
\]
It is easy to see that by monotonicity (for instance, for $j\geq n$, we have
$X_{u_{\ell}}^{2}I(|X_{u_{\ell}}|>j^{1/2})\leq X_{u_{\ell}}^{2}I(|X_{u_{\ell}%
}|>n^{1/2}))$, we have
\[
\max_{2^{r}\leq j<2^{r+1}}U_{j} \leq\frac{1}{2^{2r}}{\sum\limits_{1\leq
\ell\leq k_{2^{r+1}}}}X_{u_{\ell}}^{2}I(|X_{u_{\ell}}|>2^{r/2}) \, .
\]
Therefore, by using Markov inequality, we have to establish that
\[
\sum_{r\geq0}\frac{1}{2^{2r}}{\sum\limits_{1\leq\ell\leq k_{2^{r+1}}}%
}{\mathbb{E}} \big ( X_{u_{\ell}}^{2}I(|X_{u_{\ell}}|>2^{r/2})\big )\,<\infty
\,.
\]
or, equivalently,
\[
\sum_{n\geq1}\frac{1}{n^{3}}\sum_{\ell=1}^{k_{2n}}\mathbb{E}(X_{u_{\ell}}%
^{2}I(|X_{u_{\ell}}|>n^{1/2}))<\infty\,.
\]
This holds because of the following computation. By changing the order of
summation, and since $k_{n}\leq n^{2}$,
\begin{gather*}
\sum_{n\geq1}\frac{1}{n^{3}}{\displaystyle\sum\limits_{1\leq\ell\leq k_{2n}}%
}\mathbb{E}(X_{u_{\ell}}^{2}I(|X_{u_{\ell}}|>n^{1/2}))\leq\mathbb{E}%
\Big ({\displaystyle\sum\limits_{\ell\geq1}}X_{u_{\ell}}^{2}{\displaystyle\sum
\limits_{n\geq\sqrt{\ell}/2}}\frac{1}{n^{3}}I(|X_{u_{\ell}}|>n^{1/2})\Big )\\
\ll{\displaystyle\sum\limits_{\ell\geq1}}\frac{1}{\ell}\mathbb{E}(X_{u_{\ell}%
}^{2}I(\sqrt{2}|X_{u_{\ell}}|>\ell^{1/4}))\,.
\end{gather*}
We continue the estimate in the following way:
\[
{\displaystyle\sum\limits_{\ell\geq1}}\frac{1}{\ell}\mathbb{E}(X_{u_{\ell}%
}^{2}I(\sqrt{2}|X_{u_{\ell}}|>\ell^{1/4})) \leq{\displaystyle\sum
\limits_{\ell\geq1}}\frac{1}{\ell h(\ell^{1/4}/\sqrt{2})}\mathbb{E}%
(X_{u_{\ell}}^{2}h(|X_{u_{\ell}}|)) \ll{\displaystyle\sum\limits_{\ell\geq1}%
}\frac{1}{\ell h(\ell)}<\infty\,,
\]
where we used the fact that $h(\cdot)$ is a non-decreasing function, and
condition \eqref{higher moment}.

\smallskip

Therefore, by taking into account \eqref{aim1prth2}, to prove the theorem, it
suffices to show that
\begin{equation}
|s({\mathbf{T}}_{n})-s({\mathbf{Y}}_{n})|\rightarrow0\,\text{ a.s. as
}n\rightarrow\infty\, , \label{aim2prth2}%
\end{equation}
where $\mathbf{Y}_{n}=(Y_{u_{\ell}})_{1\leq\ell\leq k_{n}}$. With this aim, we
shall use Proposition \ref{approximation} as in the proof of Theorem
\ref{replacement}. This leads to the following estimate:
\begin{equation}
s({\mathbf{T}}_{n})-s({\mathbf{Y}}_{n}):=\mathbb{\ }R_{1,n}+\mathbb{\ }%
R_{2,n}(a)+\mathbb{\ }R_{3,n}(a)\,, \label{dec1as}%
\end{equation}
where $R_{1,n}$, $R_{2,n}(a)$ and $R_{3,n}(a)$ are respectively defined in
\eqref{defR1n}, \eqref{defR2n} and \eqref{defR3n} with the following
modifications: the $T_{n,u_{\ell}}$'s are defined by \eqref{defTul2} and the
$Z_{n,u_{\ell}}$'s are replaced by the $Y_{u_{\ell}}$'s in all the terms
involved in the decomposition.

We first prove that
\begin{equation}
|R_{1,n}|\rightarrow0\,,\text{ a.s. as $n\rightarrow\infty$}. \label{convasR1}%
\end{equation}
With this aim, as in the proof of Theorem \ref{replacement}, we use the
following decomposition:
\[
|R_{1,n}|\leq\Big |{\displaystyle\sum\limits_{1\leq\ell\leq k_{n}}}%
{\mathbb{E}}(T_{u_{\ell}}|{{\mathcal{F}}_{u_{\ell}}^{1}})\partial_{u_{\ell}%
}s(\mathbf{C}_{u_{\ell}})\Big |
+\Big |{\displaystyle\sum\limits_{1\leq\ell\leq k_{n}}}D_{u_{\ell}}%
\partial_{u_{\ell}}s(\mathbf{C}_{u_{\ell}})\Big |+\Big |{\displaystyle\sum
\limits_{1\leq\ell\leq k_{n}}}Y_{u_{\ell}}\partial_{u_{\ell}}s(\mathbf{C}%
_{u_{\ell}})\Big |\,,
\]
where $D_{u_{\ell}}:=D_{n,u_{\ell}}=T_{n,u_{\ell}}-{\mathbb{E}}(T_{n,u_{\ell}%
}|{{\mathcal{F}}_{u_{\ell}}^{1}})$. Hence, by taking into account
\eqref{boundsd} and the fact that ${\mathbb{E}}(X_{u_{\ell}}|{{\mathcal{F}%
}_{u_{\ell}}^{1}})=0$ a.s., we get that
\begin{multline}
|R_{1,n}|\leq c_{1}\frac{1}{n^{3/2}}{\displaystyle\sum\limits_{1\leq\ell\leq
k_{n}}}|{\mathbb{E}}(X_{u_{\ell}}I(|X_{u_{\ell}}|>n^{1/2})|{{\mathcal{F}%
}_{u_{\ell}}^{1}})|\\
+\Big |{\displaystyle\sum\limits_{1\leq\ell\leq k_{n}}}D_{u_{\ell}}%
\partial_{u_{\ell}}s(\mathbf{C}_{u_{\ell}})\Big |+\Big |{\displaystyle\sum
\limits_{1\leq\ell\leq k_{n}}}Y_{u_{\ell}}\partial_{u_{\ell}}s(\mathbf{C}%
_{u_{\ell}})\Big |\,. \label{decR1n*as}%
\end{multline}
We treat each term in the right hand side separately. To show that the first
term in the right-hand side converges almost surely to zero, namely:
\begin{equation}
\frac{1}{n^{3/2}}{\displaystyle\sum\limits_{1\leq\ell\leq k_{n}}}|{\mathbb{E}%
}(X_{u_{\ell}}I(|X_{u_{\ell}}|>n^{1/2})|{{\mathcal{F}}_{u_{\ell}}^{1}%
})|\rightarrow0\,\text{ a.s., as $n\rightarrow\infty$},
\label{firsttermR1nasres}%
\end{equation}
it suffices to prove (by using as before dyadic arguments), that, for any
$\varepsilon>0$,
\begin{equation}
\sum_{n\geq1}\frac{1}{n}{\mathbb{P}}\big (\max_{n\leq k<2n}\frac{1}{k^{3/2}%
}{\displaystyle\sum\limits_{1\leq\ell\leq k^{2}}}|{\mathbb{E}}(X_{u_{\ell}%
}I(|X_{u_{\ell}}|>k^{1/2})|{{\mathcal{F}}_{u_{\ell}}^{1}})|\geq\varepsilon
\big )<\infty\,. \label{firsttermR1nas}%
\end{equation}
But,
\begin{multline*}
\sum_{n\geq1}\frac{1}{n}{\mathbb{E}}\big (\max_{n\leq k<2n}\frac{1}{k^{3/2}%
}{\displaystyle\sum\limits_{1\leq\ell\leq k^{2}}}|{\mathbb{E}}(X_{u_{\ell}%
}I(|X_{u_{\ell}}|>k^{1/2})|{{\mathcal{F}}_{u_{\ell}}^{1}})|\big )\\
\leq\sum_{n\geq1}\frac{1}{n^{5/2}}{\displaystyle\sum\limits_{1\leq\ell
\leq4n^{2}}}{\mathbb{E}}(|X_{u_{\ell}}|I(|X_{u_{\ell}}|>n^{1/2}))\,,
\end{multline*}
and, since $h(\cdot)$ is a non-decreasing sequence, by (\ref{higher moment}),
\begin{multline}
\sum_{n\geq1}\frac{1}{n^{5/2}}{\displaystyle\sum\limits_{1\leq\ell\leq4n^{2}}%
}{\mathbb{E}}(|X_{u_{\ell}}|I(|X_{u_{\ell}}|>n^{1/2}))\ll{\displaystyle\sum
\limits_{\ell\geq1}}\frac{1}{\ell^{3/4}}\mathbb{E}(|X_{u_{\ell}}|I(\sqrt
{2}|X_{u_{\ell}}|>\ell^{1/4}))\label{firsttermR1nascompu}\\
\ll{\displaystyle\sum\limits_{\ell\geq1}}\frac{1}{\ell h(\ell^{1/4}/\sqrt{2}%
)}\mathbb{E}(X_{u_{\ell}}^{2}h(|X_{u_{\ell}}|))\ll{\displaystyle\sum
\limits_{\ell\geq1}}\frac{1}{\ell h(\ell)}<\infty\,.
\end{multline}
Therefore \eqref{firsttermR1nascompu} combined with Markov's inequality
implies \eqref{firsttermR1nas}, which in turn implies
\eqref{firsttermR1nasres}. We prove now that
\begin{equation}
\Big |{\displaystyle\sum\limits_{1\leq\ell\leq k_{n}}}D_{u_{\ell}}%
\partial_{u_{\ell}}s(\mathbf{C}_{u_{\ell}})\Big |\rightarrow0\,\text{ a.s. as
$n\rightarrow\infty$}. \label{secondtermR1nasres}%
\end{equation}
We start by noticing that $(D_{u_{\ell}}\partial_{u_{\ell}}s(\mathbf{C}%
_{u_{\ell}}))_{1\leq\ell\leq k_{n}}$ is a martingale difference sequence
adapted to the increasing filtration $\sigma(X_{u_{1}},\dots,X_{u_{\ell}%
},\mathbf{Y}_{n}\mathbf{)}$. Hence, by Burkholder's inequality for
complex-valued martingales (see, for instance, Lemma 2.12 Bai-Silverstein,
2010), and using (\ref{boundsd}), Cauchy-Schwartz's inequality and the
properties of conditional expectation, we obtain
\[
{\mathbb{E}}\Big |{\displaystyle\sum\limits_{1\leq\ell\leq k_{n}}}D_{u_{\ell}%
}\partial_{u_{\ell}}s(\mathbf{C}_{u_{\ell}})\Big |^{4}\ll\frac{1}%
{(n^{3/2})^{4}}\mathbb{E}\Big ({\displaystyle\sum\limits_{1\leq\ell\leq k_{n}%
}}{D}_{u_{\ell}}^{2}\Big )^{2}\ll\frac{k_{n}}{n^{6}}{\displaystyle\sum
\limits_{1\leq\ell\leq k_{n}}}\mathbb{E}({T}_{u_{\ell}}^{4})\,.
\]
By using the fact that $x^{-2}h(x)$ is non-increasing and condition
\eqref{higher moment}, we derive that
\begin{multline}
\sum_{n\geq1}\frac{k_{n}}{n^{6}}{\displaystyle\sum\limits_{1\leq\ell\leq
k_{n}}}\mathbb{E}({T}_{u_{\ell}}^{4})\ll\sum_{n\geq1}\frac{1}{n^{3}h(\sqrt
{n})}{\displaystyle\sum\limits_{1\leq\ell\leq k_{n}}}\mathbb{E}({X}_{u_{\ell}%
}^{2}h(X_{u_{\ell}}))\label{secondtermR1nasp1}\\
\ll\sum_{n\geq1}\frac{1}{nh(\sqrt{n})}\ll\sum_{n\geq1}\frac{1}{nh(n)}%
<\infty\,,
\end{multline}
which proves \eqref{secondtermR1nasres} by using Borel-Cantelli lemma. We show
now that
\begin{equation}
\Big |{\displaystyle\sum\limits_{1\leq\ell\leq k_{n}}}Y_{u_{\ell}}%
\partial_{u_{\ell}}s(\mathbf{C}_{u_{\ell}})\Big |\rightarrow0\,\text{ a.s. as
$n\rightarrow\infty$}. \label{troistermR1nasres}%
\end{equation}
To proof it we note that $(Y_{u_{\ell}}\partial_{u_{\ell}}s(\mathbf{C}%
_{u_{\ell}}))_{1\leq\ell\leq k_{n}}$ is a reversed martingale differences
sequence adapted to the decreasing filtration $\sigma(\mathbf{X}_{n}%
$,$Y_{u_{_{\ell+1}}},\dots,Y_{u_{n}})$. So, using Burkholder's inequality for
complex-valued reversed martingale differences, together with \eqref{boundsd},
we derive that
\[
{\mathbb{E}}\Big |{\displaystyle\sum\limits_{1\leq\ell\leq k_{n}}}Y_{u_{\ell}%
}\partial_{u_{\ell}}s(\mathbf{C}_{u_{\ell}})\Big |^{4}\ll\frac{1}%
{(n^{3/2})^{4}}\mathbb{E}\Big ({\displaystyle\sum\limits_{1\leq\ell\leq k_{n}%
}}{Y}_{u_{\ell}}^{2}\Big )^{2}\ll\frac{k_{n}}{n^{6}}{\displaystyle\sum
\limits_{1\leq\ell\leq k_{n}}}\mathbb{E}({Y}_{u_{\ell}}^{4})\,.
\]
But, $\mathbb{E}({Y}_{u_{\ell}}^{4})=3(\mathbb{E}({Y}_{u_{\ell}}^{2}%
))^{2}=3(\mathbb{E}({X}_{u_{\ell}}^{2}))^{2}$. Therefore
\[
\sum_{n\geq1}{\mathbb{E}}\Big |{\displaystyle\sum\limits_{1\leq\ell\leq k_{n}%
}}Y_{u_{\ell}}\partial_{u_{\ell}}s(\mathbf{C}_{u_{\ell}})\Big |^{4}\ll
\sum_{n\geq1}\frac{1}{n^{2}}<\infty\,,
\]
which proves \eqref{troistermR1nasres} by using Borel-Cantelli lemma. Starting
from \eqref{decR1n*as}, and gathering \eqref{firsttermR1nasres},
\eqref{secondtermR1nasres} and \eqref{troistermR1nasres}, the almost sure
convergence \eqref{convasR1} follows.

\smallskip

We prove now that, for any fixed positive integer $a$,
\begin{equation}
|R_{3,n}(a)|\rightarrow0\text{ a.s. as $n\rightarrow\infty$}. \label{convasR3}%
\end{equation}
By simple algebraic computations involving the inequality $b^{2}c\leq
b^{3}+c^{3}$ for any positive numbers $b$ and $c$, and the estimate of the
cardinality of $B_{u_{\ell}}(a)$ we obtain
\begin{equation}
|R_{3,n}(a)|\ll\frac{a^{2}}{n^{5/2}}{\displaystyle\sum\limits_{1\leq\ell\leq
k_{n}}}|T_{u_{\ell}}|^{3}+\frac{a^{2}}{n^{5/2}}{\displaystyle\sum
\limits_{1\leq\ell\leq k_{n}}}|Y_{u_{\ell}}|^{3}\,. \label{troisp1}%
\end{equation}
Using the fact that ${\mathbb{E}}(|Y_{u_{\ell}}|^{3})\leq2({\mathbb{E}%
}(Y_{u_{\ell}}^{2}))^{3/2}=2({\mathbb{E}}(X_{u_{\ell}}^{2}))^{3/2}$, we derive
that
\[
\sum_{n\geq1}\frac{1}{n^{7/2}}{\mathbb{E}}\Big (\max_{n\leq k<2n}%
{\displaystyle\sum\limits_{1\leq\ell\leq k^{2}}}|Y_{u_{\ell}}|^{3}\Big )
\ll\sum_{n\geq1}\frac{1}{n^{7/2}}{\displaystyle\sum\limits_{1\leq\ell
\leq4n^{2}}}{\mathbb{E}}(|Y_{u_{\ell}}|^{3})\newline\ll\sum_{n\geq1}\frac
{1}{n^{3/2}}<\infty\,,
\]
which shows, by standard arguments, that the second term in \eqref{troisp1}
converges almost surely to zero as $n\rightarrow\infty$. To end the proof of
\eqref{convasR3}, it remains to show that the first term in \eqref{troisp1}
converges almost surely to zero as $n\rightarrow\infty$. By using standard
dyadic arguments and Markov's inequality, we infer that this holds provided
that
\begin{equation}
\sum_{n\geq1}\frac{1}{n^{7/2}}{\mathbb{E}}\Big (\max_{n\leq k<2n}%
{\displaystyle\sum\limits_{1\leq\ell\leq k^{2}}}|X_{u_{\ell}}|^{3}%
I(|X_{u_{\ell}}|\leq k^{1/2})\Big )<\infty\,. \label{troisp3}%
\end{equation}
By simple computations involving the fact that $x^{-1}h(x)$ is non-increasing,
and condition \eqref{higher moment}, we get
\begin{multline}
\sum_{n\geq1}\frac{1}{n^{7/2}}{\mathbb{E}}\Big (\max_{n\leq k<2n}%
{\displaystyle\sum\limits_{1\leq\ell\leq k^{2}}}|X_{u_{\ell}}|^{3}%
I(|X_{u_{\ell}}|\leq k^{1/2})\Big )\label{troisp3**}\\
\leq\sum_{n\geq1}\frac{1}{n^{7/2}}{\displaystyle\sum\limits_{1\leq\ell
\leq4n^{2}}}{\mathbb{E}}(|X_{u_{\ell}}|^{3}I(|X_{u_{\ell}}|\leq(2n)^{1/2})
\ll\sum_{n\geq1}\frac{1}{nh(\sqrt{n})}\ll\sum_{n\geq1}\frac{1}{nh(n)}%
<\infty\,,
\end{multline}
which proves \eqref{troisp3} and ends the proof of \eqref{convasR3}.

\smallskip

It remains to handle the term $R_{2,n}(a)$ in \eqref{dec1as}. Let $\delta
\in]0,1/6[$ and, for any integer $\ell$ belonging to $[1,k_{n}]$, denote
\[
{\bar{X}}_{u_{\ell}}=X_{u_{\ell}}I(|X_{u_{\ell}}|\leq n^{\delta})\,\text{ and
}\,{\tilde{X}}_{u_{\ell}}=X_{u_{\ell}}I(n^{\delta}<|X_{u_{\ell}}|\leq
n^{1/2})\,.
\]
Using the fact that ${T}_{u_{\ell}}^{2}={\bar{X}}_{u_{\ell}}^{2}+{\tilde{X}%
}_{u_{\ell}}^{2}$, we shall use the following decomposition:%
\begin{align}
R_{2,n}(a)  &  ={\displaystyle\sum\limits_{1\leq\ell\leq k_{n}}}({\bar{X}%
}_{u_{\ell}}^{2}-\mathbb{E(}\bar{X}_{u_{\ell}}^{2}|\mathcal{F}_{u_{\ell}}%
^{a})))\partial_{u_{\ell}}^{2}s(\widetilde{\mathbf{C}}_{u_{\ell}}%
^{a})\nonumber\label{decR2as}\\
&  \quad\quad+{\displaystyle\sum\limits_{1\leq\ell\leq k_{n}}}\mathbb{(E(}%
\bar{X}_{u_{\ell}}^{2}|\mathcal{F}_{u_{\ell}}^{a})-Y_{u_{\ell}}^{2}%
)\partial_{u_{\ell}}^{2}s(\widetilde{\mathbf{C}}_{u_{\ell}}^{a}%
)+{\displaystyle\sum\limits_{1\leq\ell\leq k_{n}}}{\tilde{X}}_{u_{\ell}}%
^{2}\partial_{u_{\ell}}^{2}s(\widetilde{\mathbf{C}}_{u_{\ell}}^{a})\nonumber\\
&  :=I_{1,n}(a)+I_{2,n}(a)+I_{3,n}(a)\,.
\end{align}
By Lemma \ref{moment computation} from Section \ref{TehRes} applied with
$K=n^{\delta}$, $p=4$, $A_{u_{\ell}}=\partial_{u_{\ell}}^{2}s(\widetilde
{\mathbf{C}}_{u_{\ell}}^{a})$ (so by \eqref{boundsd}, $b_{n}= c_{2} n^{-2}$)
and ${\mathcal{G}}=\sigma({\mathbf{Y}_{n}})$, we get that
\[
\mathbb{E}|I_{1,n}(a)|^{4}\ll\frac{a^{3}}{n^{2-6\delta}}\,.
\]
Therefore, since $2-6\delta>1$,
\[
{\displaystyle\sum\limits_{n}}\mathbb{E}|I_{n}(a)|^{4} \ll a \,.
\]
So, for any positive integer $a$, the Borel-Cantelli lemma implies that
\begin{equation}
I_{1,n}(a)\rightarrow0\text{ a.s. as }n\rightarrow\infty\,. \label{asI1}%
\end{equation}
To handle the term $I_{2,n}(a)$ in \eqref{decR2as}, we apply first the
triangle inequality. Combined with \eqref{boundsd}, this leads to%
\begin{equation}
\label{termIIn}|I_{2,n}(a)|\leq\frac{1}{n^{2}}{\displaystyle\sum
\limits_{1\leq\ell\leq k_{n}}}|\mathbb{E(}\bar{X}_{u_{\ell}}^{2}%
|\mathcal{F}_{u_{\ell}}^{a})-\mathbb{E}Y_{u_{\ell}}^{2}|
+\big |{\displaystyle\sum\limits_{1\leq\ell\leq k_{n}}}\mathbb{(}Y_{u_{\ell}%
}^{2}-\mathbb{E}Y_{u_{\ell}}^{2})\partial_{u_{\ell}}^{2}s(\widetilde
{\mathbf{C}}_{u_{\ell}}^{a})\big |\,.
\end{equation}
By simple computations, we have that
\begin{multline}
\frac{1}{n^{2}}{\displaystyle\sum\limits_{1\leq\ell\leq k_{n}}}|\mathbb{E(}%
\bar{X}_{u_{\ell}}^{2}|\mathcal{F}_{u_{\ell}}^{a})-\mathbb{E}Y_{u_{\ell}}%
^{2}|\leq\frac{1}{n^{2}}{\displaystyle\sum\limits_{1\leq\ell\leq k_{n}}%
}|\mathbb{E(}X_{u_{\ell}}^{2}|\mathcal{F}_{u_{\ell}}^{a})-\mathbb{E}%
X_{u_{\ell}}^{2}|\label{I20}\\
+\frac{1}{n^{2}}{\displaystyle\sum\limits_{1\leq\ell\leq k_{n}}}%
|\mathbb{E}(X_{u_{\ell}}^{2}I(|X_{u_{\ell}}|>n^{\delta})|\mathcal{F}_{u_{\ell
}}^{a})|\,.
\end{multline}
By condition (\ref{mix a s}), the first term in (\ref{I20}) converges almost
surely to $0$ by letting first $n$ tend to infinity and then $a$ tend to
infinity. To show that the second term in (\ref{I20}) converges to zero, we
use again standard dyadic arguments and Markov's inequality, and infer that it
holds if
\begin{equation}
\sum_{n\geq1}\frac{1}{n^{3}}{\displaystyle\sum\limits_{1\leq\ell\leq4n^{2}}%
}\mathbb{E}(X_{u_{\ell}}^{2}I(|X_{u_{\ell}}|>n^{\delta}))<\infty\,.
\label{I21}%
\end{equation}
Since $h(\cdot)$ is non-decreasing, by using \eqref{higher moment}, we get
that
\begin{multline}
\sum_{n\geq1}\frac{1}{n^{3}}{\displaystyle\sum\limits_{1\leq\ell\leq4n^{2}}%
}\mathbb{E}(X_{u_{\ell}}^{2}I(|X_{u_{\ell}}|>n^{\delta}))\leq\sum_{n\geq
1}\frac{1}{n^{3}h(n^{\delta})}{\displaystyle\sum\limits_{1\leq\ell\leq4n^{2}}%
}\mathbb{E}(X_{u_{\ell}}^{2}h(|X_{u_{\ell}}|))\label{I22}\\
\leq\sum_{n\geq1}\frac{1}{nh(n^{\delta})}\leq\sum_{n\geq1}\frac{1}%
{nh(n)}<\infty\,,
\end{multline}
proving \eqref{I21}. To show that the last term in \eqref{exp} is convergent
to $0$ a.s., note that the random variables $d_{u_{\ell}}^{\prime}$ defined by
$d_{u_{\ell}}^{\prime}=\mathbb{(}Y_{u_{\ell}}^{2}-\mathbb{E}Y_{u_{\ell}}%
^{2})\partial_{u_{\ell}}^{2}s(\widetilde{\mathbf{C}}_{u_{\ell}}^{a})$ are
orthogonal. Moreover, by \eqref{boundsd}, $\mathbb{E}|d_{u_{\ell}}^{\prime
}|^{2}\ll n^{-4}\mathbb{E(}Y_{u_{\ell}}^{4})=3n^{-4}(\mathbb{E(}X_{u_{\ell}%
}^{2}))^{2}$. So
\[
{\displaystyle\sum\limits_{n\geq1}}\mathbb{E}|{\displaystyle\sum
\limits_{1\leq\ell\leq k_{n}}}d_{u_{\ell}}^{\prime}|^{2}\leq{\displaystyle\sum
\limits_{n\geq1}}{\displaystyle\sum\limits_{\ell\leq k_{n}}}\mathbb{E}%
|d_{u_{\ell}}^{\prime}|^{2}\ll{\displaystyle\sum\limits_{n\geq1}}\frac{n^{2}%
}{n^{4}}<\infty\,,
\]
which combined with the Borel-Cantelli lemma, implies that the last term in
(\ref{exp}) converges to $0$ a.s. This completes the proof of the fact that
\begin{equation}
\lim_{a\rightarrow\infty}\limsup_{n\rightarrow\infty}|I_{2,n}(a)|=0\,.
\label{asI2}%
\end{equation}
To handle the last term in \eqref{decR2as} we note that by \eqref{boundsd},
\begin{equation}
|I_{3,n}(a)|\ll\frac{1}{n^{2}}{\displaystyle\sum\limits_{1\leq\ell\leq k_{n}}%
}X_{u_{\ell}}^{2}I(|X_{u_{\ell}}|>n^{\delta})\,. \label{boundI3na1}%
\end{equation}
Using once again standard dyadic arguments and Markov's inequality, we infer
that $I_{3,n}(a)\rightarrow0\text{ a.s. as }n\rightarrow\infty$ by
\eqref{I21}. Therefore combining this fact with \eqref{asI1} and \eqref{asI2}
proves that
\begin{equation}
\lim_{a\rightarrow\infty}\limsup_{n\rightarrow\infty}|R_{2,n}(a)|=0\,.
\label{asR2n}%
\end{equation}
Finally, the decomposition \eqref{dec1as} together with \eqref{convasR1},
\eqref{convasR3} and \eqref{asR2n} implies \eqref{aim2prth2} which completes
the proof of the theorem. $\lozenge$

\subsection{Proof of Theorem \ref{stationary}}

We will follow the steps of the proof of Theorem \ref{almost sure} and in
addition we shall use the stationarity assumption and ergodic theorems. We
have to prove the counterparts of (\ref{aim1prth2}), (\ref{convasR1}),
(\ref{convasR3}), and (\ref{asR2n}). We shall just mention the differences. To
show that the almost sure convergence \eqref{aim1prth2} holds, we notice that
by taking into account \eqref{appligotze1}, it suffices to show that
\begin{equation}
\lim_{M\rightarrow\infty}\limsup_{n\rightarrow\infty}\frac{1}{n^{2}%
}{\displaystyle\sum\limits_{1\leq\ell\leq k_{n}}}X_{u_{\ell}}^{2}%
I(|X_{u_{\ell}}|>M)=0\,\text{ a.s.} \label{p1th7}%
\end{equation}
which follows by applying the ergodic theorem for stationary random fields
(see, for instance, Georgii (1988)).

Furthermore, to prove (\ref{convasR1}), we first modify the proof of
\eqref{firsttermR1nasres}. Let $M$ be a fixed positive real fixed and notice
that for any $n \geq M^{2}$,
\[
\frac{1}{n^{3/2}}{\displaystyle\sum\limits_{1\leq\ell\leq k_{n}}}|{\mathbb{E}%
}(X_{u_{\ell}}I(|X_{u_{\ell}}|>n^{1/2})|{{\mathcal{F}}_{u_{\ell}}^{1}})|
\leq\frac{1}{n^{2}}{\displaystyle\sum\limits_{1\leq\ell\leq k_{n}}}%
{\mathbb{E}}(X^{2}_{u_{\ell}}I(|X_{u_{\ell}}|>M)|{{\mathcal{F}}_{u_{\ell}}%
^{1}}) \, .
\]
Applying once again the ergodic theorem for stationary random fields, we get
\[
\lim_{M \rightarrow\infty}\limsup_{n\rightarrow\infty}\frac{1}{n^{2}%
}{\displaystyle\sum\limits_{1\leq\ell\leq k_{n}}}{\mathbb{E}}(X^{2}_{u_{\ell}%
}I(|X_{u_{\ell}}|>M)|{{\mathcal{F}}_{u_{\ell}}^{1}}) =0 \,\text{ a.s.}
\]
proving then that \eqref{firsttermR1nasres} holds. The additional change in
the proof of (\ref{convasR1}) is in the proof of \eqref{secondtermR1nasres},
and more specifically in the successive computations given in
\eqref{secondtermR1nasp1}. By taking into account the stationarity and
Fubini's theorem, we modify these computations as follows:
\[
\sum_{n\geq1}\frac{k_{n}}{n^{6}}{\displaystyle\sum\limits_{1\leq\ell\leq
k_{n}}}\mathbb{E}({T}_{u_{\ell}}^{4})\leq\mathbb{E}\Big ({X}_{{\mathbf{0}}%
}^{4}\sum_{n\geq1}\frac{1}{n^{2}}I(|X_{\mathbf{0}}|\leq n^{1/2})\Big )\ll
{\mathbb{E}}(X_{{\mathbf{0}}}^{2})<\infty\,.
\]

On another hand, to show that \eqref{convasR3} holds, the only modification
consists in the proof that the first term in the right-hand side of
\eqref{troisp1} converges almost surely to zero when $n$ to infinity. With
this aim, it suffices to write that for any positive real $M$,
\[
\frac{1}{n^{5/2}} \sum_{1 \leq\ell\leq k_{n}} |T_{u_{\ell}}|^{3} \leq
\frac{M^{3}}{n^{1/2}} + \frac{1}{n^{2}} \sum_{1 \leq\ell\leq k_{n}}%
X^{2}_{u_{\ell}}I(|X_{u_{\ell}}|>M)
\]
and to apply the ergodic theorem for stationary random fields as before
(notice that by stationarity, the second term in the right-hand side of
\eqref{troisp1} could be shown to converge almost surely to zero when $n$ to
infinity by using also the ergodic theorem).

We indicate now the differences in the proof of \eqref{asR2n}. To deal with
the first term in the right-hand side of \eqref{I20}, we notice that, by the
ergodic theorem for stationary random fields,
\[
\lim_{n\rightarrow\infty}\frac{1}{n^{2}}\sum_{1\leq j\leq i\leq n}%
|\mathbb{E}(X_{ij}^{2}-1|\mathcal{F}_{ij}^{a})|=\mathbb{E}\big (|\mathbb{E}%
(X_{\mathbf{0}}^{2}-1|{\mathcal{F}}_{\mathbf{0}}^{a})||{\mathcal{I}%
}\big )\text{ a.s.}%
\]
where ${\mathcal{I}}$ is the invariant $\sigma$-field. Note that, by
Proposition 1 in Dedecker (1998), ${\mathcal{I}}$ is included in the
${\mathbb{P}}$-completion of ${\mathcal{F}}_{\mathbf{0}}^{a}$ for all $a$.
Whence, the sequence $\mathbb{E}\big (|\mathbb{E}(X_{\mathbf{0}}%
^{2}-1|{\mathcal{F}}_{\mathbf{0}}^{a})||{\mathcal{I}}\big )_{a\geq1}$ is
almost surely decreasing, and therefore convergent almost surely. Since by
assumption, $\mathbb{E}(X_{\mathbf{0}}^{2}|{\mathcal{F}}_{\mathbf{0}}^{\infty
})=1$ a.s., by the reverse martingale theorem it follows that $\lim
_{a\rightarrow\infty}\mathbb{E}(X_{\mathbf{0}}^{2}-1|{\mathcal{F}}%
_{\mathbf{0}}^{a})=0$ a.s. and in ${\mathbb{L}}^{1}$. All these arguments
prove that the first term in the right-hand side of \eqref{I20} converges
almost surely to zero by letting first $n$ tend to infinity and after $a$ tend
to infinity. On another hand, in order to prove that the second term in the
right-hand side of \eqref{I20} converges almost surely to zero when $n$ tends
to infinity, it suffices to show that, for any positive integer $a$,
\[
\lim_{M\rightarrow\infty}\limsup_{n\rightarrow\infty}\frac{1}{n^{2}%
}{\displaystyle\sum\limits_{1\leq\ell\leq k_{n}}}|\mathbb{E}(X_{u_{\ell}}%
^{2}I(|X_{u_{\ell}}|>M)|{\mathcal{F}}_{u_{\ell}}^{a})|=0\ \text{a.s.},
\]
which follows by the ergodic theorem for stationary random fields. Similarly,
the ergodic theorem for stationary random fields together with the bound in
\eqref{boundI3na1} allows us to prove that $I_{3,n}(a)$ converges almost
surely to zero when $n$ tends to infinity.

Therefore, under the conditions of Theorem \ref{stationary}, the conclusion of
Theorem \ref{almost sure} holds. Furthermore, condition (\ref{Bogus0}) is
satisfied, hence the result follows from Corollary \ref{semicircular}.
$\lozenge$

\section{Technical Results\label{TehRes}}

Below we give an approximation theorem needed for the proof of the main
theorems. A related approximation result is in Chatterjee (2006).

\begin{proposition}
\label{approximation}Suppose that $\mathbf{X}:=(X_{1},\dots,X_{m})$ and
$\mathbf{Z}:=(Z_{1},\dots,Z_{m})$ are random vectors in ${\mathbb{R}}^{m}$.
Suppose that $f:{\mathbb{R}}^{m}\rightarrow{\mathbb{C}}$ is a function three
times differentiable with bounded partial derivatives
\[
|\partial_{uvw}^{3}f(x)|\leq L_{3}\text{ \ for all }x\text{ all }u,\text{
}v,\text{ }w.
\]
Let $B_{k}$ be a subset of the set $\{i\in{\mathbb{N}}\,:\,1\leq i<k\}$.
Denote by $B_{k}^{c}=\{i\in{\mathbb{N}}\,:\,1\leq i<k\}{\small \setminus
B_{k}}$. Define a vector $\mathbf{U}_{k-1}=\mathbf{(}U_{1},\dots,U_{k-1})$
such that $U_{i}=0$ if $i\in B_{k}$ and $U_{i}=X_{i}$ if $i\in B_{k}^{c}$.
Then%
\[
f(\mathbf{X})-f(\mathbf{Z})=R_{1}+R_{2}+R_{3}%
\]
where%
\[
R_{1}={\displaystyle\sum\limits_{1\leq k\leq m}}\mathbb{(}X_{k}-Z_{k}%
)\partial_{k}f(X_{1},\dots,X_{k-1},0,Z_{k+1},\dots,Z_{m}),
\]
\[
R_{2}=\frac{1}{2}{\displaystyle\sum\limits_{1\leq k\leq m}}\mathbb{(}X_{k}%
^{2}-Z_{k}^{2})\partial_{k}^{2}f(\mathbf{U}_{k-1},0,Z_{k+1},\dots,Z_{m}),
\]
and%
\[
|R_{3}|\leq L_{3}{\displaystyle\sum\limits_{1\leq k\leq m}}\mathbb{(}X_{k}%
^{2}+Z_{k}^{2}){\displaystyle\sum\limits_{u\in B_{k}}}|X_{u}|+L_{3}%
{\displaystyle\sum\limits_{1\leq k\leq m}}|X_{k}|^{3}+L_{3}{\displaystyle\sum
\limits_{1\leq k\leq m}}|Z_{k}|^{3}.
\]

\end{proposition}

\noindent\textbf{Proof.} For any $k\in\{0,\dots,m\}$, we define the following
vectors
\[
Y_{k}=(X_{1},\dots,X_{k},Z_{k+1},\dots,Z_{m})\,\text{ and }\,Y_{k}%
^{(0)}=(X_{1},\dots,X_{k-1},0,Z_{k+1},\dots,Z_{m})\,.
\]
Then, we have the telescoping decomposition:
\[
f(\mathbf{X})-f(\mathbf{Z})={\displaystyle\sum\limits_{1\leq k\leq m}}%
(f(Y_{k})-f(Y_{k}^{0})+f(Y_{k}^{0})-f(Y_{k-1}))\,.
\]
By applying the Taylor expansion of order two, we get
\[
f(Y_{k})-f(Y_{k}^{0})=X_{k}\partial_{k}f(Y_{k}^{0})+\frac{1}{2}X_{k}%
^{2}\partial_{k}^{2}f(Y_{k}^{0})+R_{3}^{\prime}\,,
\]
where $|R_{3}^{\prime}|\leq L_{3}|X_{k}|^{3}$. By writing a similar expansion
for $f(Y_{k-1})-f(Y_{k}^{0})$ leads to
\begin{equation}
f(\mathbf{X})-f(\mathbf{Z})={\displaystyle\sum\limits_{1\leq k\leq m}%
}\mathbb{(}X_{k}-Z_{k})\partial_{k}f(Y_{k}^{0})+\frac{1}{2}{\displaystyle\sum
\limits_{1\leq k\leq m}}\mathbb{[}X_{k}^{2}-Z_{k}^{2}]\partial_{k}^{2}%
f(Y_{k}^{0})+R_{3}^{\prime\prime}\,, \label{expansion}%
\end{equation}
where%
\[
|R_{3}^{\prime\prime}|\leq L_{3}{\displaystyle\sum\limits_{1\leq k\leq m}%
}(|X_{k}|^{3}+|Z_{k}|^{3})\,.
\]
We continue to estimate the second term in the right-hand side of
(\ref{expansion}). Let $V_{k}=X_{k}^{2}-Z_{k}^{2}$ and write
\begin{multline*}
V_{k}\partial_{k}^{2}f(Y_{k}^{(0)})=V_{k}\partial_{k}^{2}f(\mathbf{U}%
_{k-1},0,Z_{k+1},\dots,Z_{m})\\
+(V_{k}\partial_{k}^{2}f(Y_{k}^{(0)})-V_{k}\partial_{k}^{2}f(\mathbf{U}%
_{k-1},0,Z_{k+1},\dots,Z_{m}))\,.
\end{multline*}
By Taylor expansion of first order and taking into account the bounds for the
derivatives, we have
\[
|V_{k}\partial_{k}^{2}f(Y_{k}^{(0)})-V_{k}\partial_{k}^{2}f(\mathbf{U}%
_{k-1},0,Z_{k+1},\dots,Z_{m})|\leq L_{3}|V_{k}|\sum_{u\in B_{k}}|X_{u}|.
\]
Finally set
\[
R_{3}=R_{3}^{\prime\prime}+\frac{1}{2}{\displaystyle\sum\limits_{1\leq k\leq
m}}(V_{k}\partial_{k}^{2}f(Y_{k}^{(0)})-V_{k}\partial_{k}^{2}f(\mathbf{U}%
_{k-1},0,Z_{k+1},\dots,Z_{m}))\,,
\]
and the result follows. $\lozenge$

\bigskip

We state next Lemma 2.1 in G\"{o}tze \textit{et al.} (2012).

\begin{lemma}
\label{lmagotze} Let $\mathbf{x}=(x_{ij})_{1\leq j\leq i\leq n}$ and
$\mathbf{y}=(y_{ij})_{1\leq j\leq i\leq n}$ two elements of ${\mathbb{R}%
}^{k_{n}}$ where $k_{n}=n(n+1)/2$. Let $z=u+iv\in\mathbb{C}^{+}$ and $s (
\cdot) := s ( \cdot, z)$ be the function from $\mathbb{R}^{k_{n}}$ to
$\mathbb{C}$ defined by \eqref{defff}. Then
\[
|s(\mathbf{x})-s(\mathbf{y})|\leq\frac{1}{v^{2}} \big (\frac{1}{n^{2}}%
\sum_{i=1}^{n}(x_{ii}-y_{ii})^{2}+2\sum_{i=1}^{n}\sum_{j=1}^{i-1}%
(x_{ij}-y_{ij})^{2} \big )^{1/2}\,.
\]

\end{lemma}

The following lemma is an easy consequence of the well-known Gaussian
interpolation. For reference we cite Talagrand (2010) Section 1.3, Lemma 1.3.1.

\begin{lemma}
\label{inter}Suppose that $\mathbf{Y}=(Y_{1}, \dots, Y_{m})$ and
$\mathbf{Z}=(Z_{1}, \dots, Z_{m})$ are Gaussian centered random vectors in
${\mathbb{R}}^{m}$ with independent components. Suppose that $f:{\mathbb{R}%
}^{m}\rightarrow{\mathbb{C}}$ is a function twice differentiable with bounded
partial derivatives%
\[
\text{\ }|\partial_{u}f(\mathbf{x})|\leq L_{1}\text{\ and }|\partial_{u}%
^{2}f(\mathbf{x})|\leq L_{2}\text{ for all }\mathbf{x},\text{ }u.
\]
Then%
\[
|\mathbb{E}f(\mathbf{Y)-}\mathbb{E}f(\mathbf{Z)|}\leq\frac{L_{2}}{2}\sum
_{i=1}^{n} |\mathbb{E}Y_{i}^{2}-\mathbb{E}Z_{i}^{2}| \, .
\]

\end{lemma}

In the next lemma we compute moments of some terms which appear in the proofs
of Theorems \ref{replacement} and \ref{almost sure}. Before stating it, for
reader convenience, let us recall some notations: $k_{n}=n(n+1)/2$ and
$(u_{\ell},1\leq\ell\leq k_{n})$ are double indexes ordered in the strict
lexicographic order. To be more precise, for any integer $\ell\in
\lbrack1,k_{n}]$, if $i$ is the integer in $[1,n]$ such that $\frac{i(i-1)}%
{2}+1\leq\ell\leq\frac{i(i+1)}{2}$, then $\ell=\frac{i(i-1)}{2}+j$ with
$j\in\{1,\dots,i\}$ and $u_{\ell}=(i,j)$.

\begin{lemma}
\label{moment computation} Let $a$ and $K$ be two positive integers. For any
integer $\ell\in\lbrack1,k_{n}]$, let
\[
{\bar{X}}_{u_{\ell}}=X_{u_{\ell}}I(|X_{u_{\ell}}|\leq K)\,.
\]
Let ${\mathcal{G}}$ be a sigma algebra independent of $\sigma\{(X_{ij}%
\,)_{i,j\in\mathbb{Z}^{2}}\}$ and $\mathcal{F}_{u_{\ell}}^{a}$ be defined by
\eqref{deffiltration}. Let $(A_{u_{\ell}})_{1\leq\ell\leq k_{n}}$ be a
sequence of complex-valued random variables such that $A_{u_{\ell}}$ is
$\mathcal{F}_{u_{\ell}}^{a}\vee{\mathcal{G}}$-measurable and
\[
\max_{1\leq\ell\leq k_{n}}|A_{u_{\ell}}|\leq b_{n}\,\text{\ a.s.}%
\]
Assume that condition \eqref{bound} holds. Then for any $p\geq2$,
\[
\mathbb{E}\big |{\sum\limits_{1\leq\ell\leq k_{n}}}\mathbb{(}{\bar{X}%
}_{u_{\ell}}^{2}-\mathbb{E(}{\bar{X}}_{u_{\ell}}^{2}|\mathcal{F}_{u_{\ell}%
}^{a}))A_{u_{\ell}})\big |^{p}\ll K^{2(p-1)}b_{n}^{p}(a^{p/2}n^{3p/2}%
+a^{p-1}n^{p+1})\,.
\]

\end{lemma}

\noindent\textbf{Proof.} The proof is based Burkholder's inequality for
differences of martingale with complex valued random variables. Because the
filtration $\mathcal{F}_{u_{\ell}}^{a}$ is not nested we shall apply a
blocking procedure. Let $v_{n}=[n/a]$ where $[x]$ denotes the integer part of
$x$. Setting
\[
d_{i,j}=\mathbb{(}{\bar{X}}_{ij}^{2}-\mathbb{E(}{\bar{X}}_{ij}^{2}%
|\mathcal{F}_{ij}^{a}))A_{ij}\,,\text{ if }1\leq j\leq i\leq n
\]
and
\[
d_{i,j}=0\,,\text{ if }1\leq i<j\leq n.
\]
For pointing out an adapted martingale structure, we decompose the sum in the
following way:
\[
{\sum\limits_{1\leq\ell\leq k_{n}}}\mathbb{(}{\bar{X}}_{u_{\ell}}%
^{2}-\mathbb{E(}{\bar{X}}_{u_{\ell}}^{2}|\mathcal{F}_{u_{\ell}}^{a}%
))A_{u_{\ell}}=\sum_{m=1}^{a}\sum_{k=0}^{v_{n}-1}\sum_{j=1}^{ka+m}%
d_{ka+m,j}+\sum_{i=v_{n}a+1}^{n}\sum_{j=1}^{i}d_{i,j}\,,
\]
implying that
\begin{equation}
\Big \Vert{\sum\limits_{1\leq\ell\leq k_{n}}}\mathbb{(}{\bar{X}}_{u_{\ell}%
}^{2}-\mathbb{E(}{\bar{X}}_{u_{\ell}}^{2}|\mathcal{F}_{u_{\ell}}%
^{a}))A_{u_{\ell}}\Big \Vert_{p}\leq\sum_{m=1}^{a}\sum_{j=1}^{n}%
\Big \Vert\sum_{k=0}^{v_{n}-1}d_{ka+m,j}\Big \Vert_{p}+\sum_{i=v_{n}a+1}%
^{n}\sum_{j=1}^{i}\Vert d_{i,j}\Vert_{p}\,. \label{p1lmacompmom}%
\end{equation}
To handle the first term in the right-hand side of the above inequality, we
note that for $m$ and $j$ fixed, $(d_{ka+m,j})_{k\geq0}$ is a complex-valued
sequence of martingale differences with respect to the filtration
$\mathcal{F}_{ka+m,j}^{0}\vee{\mathcal{G}}$. To see this, just note that
$d_{ka+m,j}$ is adapted to $\mathcal{F}_{ka+m,j}^{0}\vee{\mathcal{G}}$ and we
also have, for $k\geq1$, $\mathcal{F}_{(k-1)a+m,j}^{0}\subset\mathcal{F}%
_{ka+m,j}^{a}$. Then, using also that $A_{ka+m,j}$ is $\mathcal{F}%
_{ka+m,j}^{a}\vee{\mathcal{G}}$-measurable and that ${\mathcal{G}}$ is
independent of $\sigma(X_{u_{i}}\,,\,1\leq i\leq k_{n})$, we get for $k\geq
1$,
\begin{multline*}
\mathbb{E}(d_{ka+m,j}|\mathcal{F}_{(k-1)a+m,j}^{0}\vee{\mathcal{G}})=\\
\mathbb{E}(A_{ka+m,j}\mathbb{E(}X_{ka+m,j}^{2}-\mathbb{E}(X_{ka+m,j}%
^{2}|\mathcal{F}_{ka+m,j}^{a})|\mathcal{F}_{ka+m,j}^{a})|\mathcal{F}%
_{(k-1)a+m,j}^{0}\vee{\mathcal{G}})=0\text{ a.s.}%
\end{multline*}
Therefore, by applying Burkholder's inequality for differences of martingale
with complex valued (see, for instance, Lemma 2.12 in Bai-Silverstein, 2010),
it follows that there exists a universal positive constant $C_{p}$ depending
only on $p$ such that, for any $m\in\{1,\dots,a\}$ and any $j\in\{1,\dots
,n\}$,
\begin{equation}
\Big \Vert\sum_{k=0}^{v_{n}-1}d_{ka+m,j}\Big \Vert_{p}^{p}\leq C_{p}%
\Big \Vert\sum_{k=0}^{v_{n}-1}|d_{ka+m,j}|^{2}\Big \Vert_{p/2}^{p/2}\,.
\label{p3lmacompmom}%
\end{equation}
But, for any $1\leq j\leq i\leq n$,
\begin{equation}
|d_{i,j}|\leq b_{n}|{\bar{X}}_{ij}^{2}-\mathbb{E(}{\bar{X}}_{ij}%
^{2}|\mathcal{F}_{ij}^{a})|\leq2b_{n}K^{2}\,, \label{bounddkj}%
\end{equation}
implying that
\[
\Big \vert\sum_{k=0}^{v_{n}-1}|d_{ka+m,j}|^{2}\Big \vert^{p/2}\leq
2^{p-1}K^{2(p-1)}b_{n}^{p-1}\sum_{k=0}^{v_{n}-1}|d_{ka+m,j}|\,.
\]
Hence, starting from \eqref{p3lmacompmom} and using the above upper bound, we
get
\[
\Big \Vert\sum_{k=0}^{v_{n}-1}d_{ka+m,j}\Big \Vert_{p}^{p}\leq C_{p}%
2^{p-1}K^{2(p-1)}v_{n}^{(p-2)/2}b_{n}^{p-1}\sum_{k=0}^{v_{n}-1}{\mathbb{E}%
}(|d_{ka+m,j}|)\,.
\]
which combined with the first part of \eqref{bounddkj} entails
\[
\Big \Vert\sum_{k=0}^{v_{n}-1}d_{ka+m,j}\Big \Vert_{p}^{p}\leq C_{p}%
2^{p}K^{2(p-1)}b_{n}^{p}v_{n}^{(p-2)/2}\sum_{k=0}^{v_{n}-1}{\mathbb{E}%
}(X_{ka+m,j}^{2}){\mathbf{1}}_{ka+m\geq j}\,.
\]
Therefore, using H\"{o}lder's inequality and the above inequality, we derive
\begin{align*}
\Big (\sum_{m=1}^{a}\sum_{j=1}^{n}\Big \Vert\sum_{k=0}^{v_{n}-1}%
d_{ka+m,j}\Big \Vert_{p}\Big )^{p}  &  \leq(an)^{p-1}C_{p}2^{p}K^{2(p-1)}%
b_{n}^{p}v_{n}^{(p-2)/2}\sum_{m=1}^{a}\sum_{k=0}^{v_{n}-1}\sum_{j=1}%
^{ka+m}{\mathbb{E}}(X_{ka+m,j}^{2})\\
&  \leq(an)^{p-1}C_{p}2^{p}K^{2(p-1)}b_{n}^{p}v_{n}^{(p-2)/2}\sum_{i=1}%
^{n}\sum_{j=1}^{i}{\mathbb{E}}(X_{ij}^{2})\,.
\end{align*}
By taking into account condition (\ref{bound}) and the fact that $v_{n}\leq
n/a$, it follows that
\begin{equation}
\Big (\sum_{m=1}^{a}\sum_{j=1}^{n}\Big \Vert\sum_{k=0}^{v_{n}-1}%
d_{ka+m,j}{\mathbf{1}}_{ka+m\geq j}\Big \Vert_{p}\Big )^{p}\ll a^{p/2}%
K^{2(p-1)}b_{n}^{p}n^{3p/2}\,. \label{p4lmacompmom}%
\end{equation}
We handle now the second term in the right-hand side of inequality
\eqref{p1lmacompmom}. With this aim, we use H\"{o}lder's inequality and
\eqref{bounddkj} to get
\begin{align*}
\Big (\sum_{i=v_{n}a+1}^{n}\sum_{j=1}^{i}\Vert d_{i,j}\Vert_{p}\Big)^{p}  &
\leq(n-v_{n}a)^{p-1}n^{p-1}\sum_{i=v_{n}a+1}^{n}\sum_{j=1}^{i}\mathbb{E}%
(|d_{i,j}|^{p})\\
&  \leq a^{p-1}n^{p-1}(2b_{n}K^{2})^{p-1}(2b_{n})\sum_{i=1}^{n}\sum_{j=1}%
^{i}\mathbb{E}(X_{ij}^{2})\,.
\end{align*}
Hence condition (\ref{bound}) implies
\begin{equation}
\Big (\sum_{i=v_{n}a+1}^{n}\sum_{j=1}^{i}\Vert d_{i,j}\Vert_{p}\Big)^{p}\ll
a^{p-1}K^{2(p-1)}b_{n}^{p}n^{p+1}\,. \label{p5lmacompmom}%
\end{equation}
The lemma follows by taking into account the upper bounds \eqref{p4lmacompmom}
and \eqref{p5lmacompmom} in \eqref{p1lmacompmom}. $\lozenge$

\begin{remark}
\label{remark4} Our proof shows that the conclusion of the lemma still holds
if we replace the filtration $\mathcal{F}_{ij}^{a}$ by a larger filtration
$\mathcal{K}_{ij}^{a}$ (for $a \geq0$ fixed) with the following properties:
for any $j \leq i$, $\mathcal{F}_{ij}^{0}\subseteq\mathcal{K}_{ij}^{0}$,
$\mathcal{K}_{ij}^{a}\subseteq\mathcal{K}_{i j}^{0}$, $\mathcal{K}_{ij}%
^{0}\subseteq\mathcal{K}_{i+1, j}^{0}$ and $\mathcal{K}_{i-a,j}^{0}%
\subseteq\mathcal{K}_{ij}^{a}$ for $i \geq a+1$. Moreover in the statement of
the lemma, the filtration ${\mathcal{G}}$ has to be assumed to be independent
of $\sigma( \bigcup_{i,j}\mathcal{K}_{ij}^{0})$. For instance, we can take
${\mathcal{K}}_{ij}^{a}=\sigma(X_{uv}\,:\,(u,v)\in B_{ij}^{a})$ where
$B_{ij}^{a}$ is defined in \eqref{defB}.
\end{remark}


\noindent

\begin{thebibliography}{99}                                                                                               %


\bibitem {A}Adamczak, R. (2011). On the Marchenko-Pastur and circular laws for
some classes of random matrices with dependent entries. \textit{Electronic
Journal of Probability} \textbf{16 }1065-1095.

\bibitem {A2}Adamczak, R. (2013). Some remarks on the Dozier-Silverstein
theorem for random matrices with dependent entries. \textit{Random Matrices:
Theory. Appl.} \textbf{02} 1-46.

\bibitem {a}Arnold, L. (1971). On Wigner's semicircle law for the eigenvalues
of random matrices. \textit{Z.Wahrscheinlichkeitstheorie und Verw. Gebiete}
\textbf{19} 191-198.

\bibitem {BS}Bai, Z. and J.W. Silverstein (2010). \textit{Spectral analysis of
large dimensional random matrices}. Springer, New York, second edition.

\bibitem {BM}Boutet de Monvel A. and A. Khorunzhy (1999). On the norm and
eigenvalue distribution of large random matrices. \textit{Ann. Probab.}
\textbf{27} 913-944.

\bibitem {Ct}Chatterjee, S. (2006). A generalization of the Lindeberg
principle.\textit{ Ann. Probab.} \textbf{34} 2061-2076.

\bibitem {De98}Dedecker, J. (1998). A central limit theorem for stationary
random fields.\textit{ Probab. Theory Related Fields} \textbf{110} 397-426.

\bibitem {DT}Doukhan, P. and L. Truquet. (2007). A fixed point approach to
model random fields. \textit{Alea} \textbf{3} 111-132.

\bibitem {Ge}Georgii, H.O. (1988). \textit{Gibbs Measures and Phase
Transitions}. De Gruyter.

\bibitem {GH}Geronimo J.S. and T. Hill. (2003). Necessary and Sufficient
Condition that the Limit of Stieltjes Transforms is a Stieltjes Transform.
\textit{J. of Approx. Theory} \textbf{121} 54-60.

\bibitem {gg}Girko, V. (2013). The Generalized Circular Law. \textit{Random
Operators and Stochastic Equations} \textbf{21} 67--109.

\bibitem {g}Girko, V., Kirsch, W. and A. Kutzelnigg (1994). A necessary and
sufficient condition for the semicircle law. \textit{Random Operators and
Stochastic Equations} \textbf{2} 195-202.

\bibitem {GT2}G\"{o}tze, F. and A. Tikhomirov (2004). Limit theorems for
spectra of positive random matrices under dependence. \textit{Zap. Nauchn.
Sem.} S.-Peterburg. Otdel. Mat. Inst. Steklov. (POMI),\textbf{ 311} (Veroyatn.
i Stat. 7) 92-23, 299.

\bibitem {GT}G\"{o}tze, F. and A.N. Tikhomirov (2006). Limit theorems for
spectra of random matrices with martingale structure. \textit{Teor. Veroyatn.
Primen}., \textbf{51} 171-192.

\bibitem {GNT}G\"{o}tze, F., Naumov, A. and A. Tikhomirov (2012). Semicircle
law for a class of random matrixes with dependent entries. arXiv:math/0702386v1

\bibitem {MP}Marchenko, V.A. and L.A. Pastur. (1967). Distribution for some
sets of random matrices. \textit{Math. USSR-Sb}. \textbf{1} 457--483.

\bibitem {r}O'Rourke, S. (2012). A note on the Marchenko-Pastur law for a
class of random matrices with dependent entries. arXiv:1201.3554.

\bibitem {P}Pastur, L.A. (1973). Spectra of random selfadjoint
operators.\textit{ Uspehi Mat. Nauk} \textbf{28} 3-64.

\bibitem {Tal}Talagrand M. (2010). \textit{Mean Field Models for Spin Glasses.
Vol 1. Basic Examples.} Springer.

\bibitem {T}Tao, T and V. Vu (2010). Random matrices: universality of ESDs and
the circular law. With an appendix by Manjunath Krishnapur. \textit{Ann.
Probab}. \textbf{38} 2023--2065.

\bibitem {Wa}Wachter, K.W. (1976). Proceedings of the Computer Science and
Statistics 9th Annual Symposium on the Interface, pp. 299--308.

\bibitem {W}Wigner E.P. (1958). On the distribution of the roots of certain
symmetric matrices.\textit{ Ann. of Math.} \textbf{67} 325-327.
\end{thebibliography}
\end{document}